\documentclass[aap,seceqn,nameyear,MSNbibl,dvips]{arximspdf}

\doi{10.1214/09-AAP674}
\volume{20}
\issue{5}
\pubyear{2010}
\firstpage{1854}
\lastpage{1890}

\makeatletter

\newcommand{\PH}{\mathit{Ph}}
\newcommand{\GI}{\mathit{GI}}

\newtheorem{lemma}{Lemma}
\newtheorem{theorem}{Theorem}
\newproclaim{definition}{Definition}
\newproclaim{remark}{Remark}

\makeatother

\begin{document}
\begin{frontmatter}

\title{Many-server diffusion limits for $G/\mathit{P\lowercase{h}}/\lowercase{n}+\GI$ queues\thanksref{T1}}
\runtitle{Many-server diffusion limits}

\thankstext{T1}{Supported in part by NSF Grants CMMI-0727400 and
CMMI-0825840 and by an IBM Faculty Award.}

\begin{aug}
\author[A]{\fnms{J. G.} \snm{Dai}\corref{}\ead[label=e1]{dai@gatech.edu}},
\author[A]{\fnms{Shuangchi} \snm{He}\ead[label=e2]{heshuangchi@gatech.edu}} and
\author[B]{\fnms{Tolga} \snm{Tezcan}\ead[label=e3]{ttezcan@uiuc.edu}\ead[label=e4]{Tolga.Tezcan@simon.rochester.edu}}
\runauthor{J. G. Dai, S. He and T. Tezcan}
\affiliation{Georgia Institute of Technology, Georgia Institute of
Technology\break and University
of Illinois at Urbana-Champaign}
\address[A]{J. G. Dai\\
S. He\\
H. Milton Stewart School of Industrial\\
\quad and Systems Engineering\\
Georgia Institute of Technology\\
Atlanta, Georgia 30332\\
USA\\
\printead{e1}\\
\phantom{E-mail: }\printead*{e2}} 
\address[B]{T. Tezcan\\
Department of Industrial\\
\quad and Enterprise Systems Engineering\\
University of Illinois at Urbana-Champaign\\
Urbana, Illinois 61801\\
USA\\
\printead{e3}\\
and\\
Simon Graduate School of Business\\
University of Rochester\\
305 Schlegel Hall\\
Rochester, New York 14627\\
USA\\
\printead{e4}}
\end{aug}

\received{\smonth{12} \syear{2008}}
\revised{\smonth{11} \syear{2009}}

%
\begin{abstract}
This paper studies many-server limits for multi-server queues that
have a phase-type service time distribution and allow for customer
abandonment. The first set of limit theorems is for critically
loaded $G/\PH/n+\GI$ queues, where the patience times are independent and
identically distributed following a general distribution. The next
limit theorem is for overloaded $G/\PH/n+M$ queues, where the patience
time distribution is restricted to be exponential. We prove that a
pair of diffusion-scaled total-customer-count and server-allocation
processes, properly centered, converges in distribution to a
continuous Markov process as the number of servers $n$ goes to
infinity. In the overloaded case, the limit is a multi-dimensional
diffusion process, and in the critically loaded case, the limit is a
simple transformation of a diffusion process. When the queues are
critically loaded, our diffusion limit generalizes the result by
\citet{PuhalskiiReiman00} for $\GI/\PH/n$ queues without customer
abandonment. When the queues are overloaded, the diffusion limit
provides a refinement to a fluid limit and it generalizes a result by
\citet{Whitt04a} for $M/M/n/+M$ queues with an exponential
service time
distribution. The proof techniques employed in this paper are
innovative. First, a perturbed system is shown to be equivalent to
the original system. Next, two maps are employed in both fluid and
diffusion scalings. These maps allow one to prove the limit theorems
by applying the standard continuous-mapping theorem and the standard
random-time-change theorem.
\end{abstract}

%

\end{frontmatter}

\section{Introduction}
\label{sec:intro}

This paper studies many-server limits for multi-server queues that
allow for
customer abandonment. These queues are assumed to have a phase-type service
time distribution. We consider two separate parameter regimes: one for
critically loaded many-server queues and the other for overloaded
many-server queues.

As argued in the seminal paper of \citet{HalWhi81}, for a
critically loaded
many-server queue, the system provides high-quality service and at the same
time achieves high server utilization. Thus, the critically loaded parameter
regime is also known as the Quality and Efficiency-Driven (QED) limiting
regime or the Halfin--Whitt limiting regime. For the overloaded $M/M/n+M$
model, \citet{Whitt04a} demonstrates that a certain fluid
approximation can be
useful in predicting the \textit{steady-state} performance of the multi-server
system. He further demonstrates that a diffusion limit provides a
refined approximation.

Our first set of results is for critically loaded $G/\PH/n+\GI$ queues, whose
patience times are independent and identically distributed (i.i.d.)
following a
general distribution. In Theorem \ref{thm:critical}, we prove that a
pair of
diffusion-scaled \textit{total-customer-count} and \textit{server-allocation
processes} converges in distribution to a continuous Markov process
$(\tilde{X},\tilde{Z})$. In Theorem \ref{cor:1}, we prove that the
diffusion-scaled
\textit{customer-count-vector process} converges to a diffusion process
$\tilde{Y}$. In \mbox{Theorem~\ref{thm:virtualwaiting}}, the diffusion-scaled
\textit{virtual
waiting time process} converges in distribution to a constant multiple of
$(\tilde{X})^{+}$, which serves as the limit of the diffusion-scaled
queue-length process. Our second result is for overloaded $G/\PH/n+M$
queues, whose patience time distribution is restricted to be
exponential. In
Theorem \ref{thm:overloaded}, we prove that the pair of diffusion-scaled
\textit{total-customer-count} and \textit{server-allocation processes} converges
in distribution to a diffusion process. Although the limit $(\tilde{X}%
,\tilde{Z})$ in Theorem \ref{thm:critical} is not a diffusion process
in a
strict sense (see discussions below the statement of Theorem \ref
{cor:1}), we
still call it a diffusion limit because it is a simple transformation
of the
diffusion process $\tilde{Y}$ in Theorem \ref{cor:1}. This
terminology is
consistent with the usage in conventional heavy traffic, where the limit
process is often a constrained diffusion process [see, e.g.,
\citet{rei84}].

In the critically-loaded regime, \citet{HalWhi81} is the first
paper to establish a
diffusion limit for the $\GI/M/n$ model. \citet{PuhalskiiReiman00}
establish a
diffusion limit for the $\GI/\PH/n$ model, where the service time distribution
is phase-type. \citet{GMR02} prove a diffusion limit for the
$M/M/n+M$ model,
which allows for customer abandonment. \citet{Whitt05a}
generalizes the result to
the $G/M/n+M$ model. In the same paper, Whitt proves a stochastic-process
limit for the $G/H_{2}^{\ast}/n$ model; this limit is not a diffusion process
although a simple transformation of it is a diffusion process. Our
first set
of results, Theorems \ref{thm:critical} and \ref{cor:1}, extends the
result of
\citet{PuhalskiiReiman00} to the $G/\PH/n+\GI$ model, which allows
for customer
abandonment with a general patience time distribution. It also extends the
work of \citet{GMR02} and \citet{Whitt05a} to allow for
phase-type service time
distributions. For the overloaded $G/\PH/n+M$ model, Theorem
\ref{thm:overloaded} generalizes \citet{Whitt05a} to the
$G/\PH/n+M$ model. The
diffusion limit in the theorem provides a refinement to a fluid limit.

In addition to these limit theorems, the techniques used in the proofs are
innovative. When the patience time distribution is exponential, we first
establish a sample-path representation for our $G/\PH/n+M$ model. This
representation allows us to obtain the total-customer-count and the
server-allocation processes as a map of primitive processes with a
random time change. These primitive processes are either assumed or
known to
satisfy functional central limit theorems (FCLTs). Therefore, our limits
follow from the standard continuous-mapping theorem and the standard
random-time-change theorem [see, e.g., \citet{ethkur86} and
\citet{bil99}]. When the queues are critically loaded, a result of
\citet{DaiHe09} proves that the performance of these queues is not
sensitive to
the distribution of patience times, thus allowing us to prove
Theorems \ref{thm:critical}--\ref{thm:virtualwaiting} for general patience
time distributions with a mild regularity condition.

\citet{HalWhi81}, \citet{GMR02} and \citet{Whitt04a} all
use Stone's theorem to
prove diffusion limit theorems in the critically-loaded regime.
\citet{Stone63} is set
up for convergence of Markov chains to a diffusion process. This
setting makes
the generalization to nonrenewal arrival processes difficult.
\citet{PuhalskiiReiman00} also use the continuous-mapping approach
for the
$\GI/\PH/n$ model. They employ a different sample-path representation for the
total-customer-count process. Their representation requires them to use
extensively \textit{martingale FCLTs} in their proofs, whereas our
approach uses
standard FCLTs for random walks and Poisson processes. \citet
{Whitt07c} review
a number of sample-path representations and martingale proofs for many-server
heavy traffic limits, and \citet{Whitt07b} surveys the proof
techniques for
establishing martingale FCLTs. Our proofs show that for multi-server queues
with a phase-type service time distribution and an exponential patience time
distribution, there is a general approach to proving limit theorems, without
employing martingale FCLTs. Note that when the patience time
distribution is general, our proofs for the diffusion limits in the
critically-loaded
regime rely on a result of \citet{DaiHe09}, which is proved by using
a martingale FCLT.

Our sample-path representation is based on the equivalence of our multi-server
system to a perturbed system as illustrated in \citet{Tezcan06}. This
representation has been used in \citet{DaiTezcan05}, \citet
{DaiTez06} and
\citet{daiTez07a} for multi-server-pool systems when service and
patience times
have exponential distributions. The sample-path argument has been explored
previously in the setting of Markovian networks in \citet
{MandelStrong} for
strong approximations and \citet{MandelbaumPats98} for general
state-dependent networks.

In our continuous-mapping approach, we have heavily exploited some maps from
$\mathbb{D}^{K+1}$ to $\mathbb{D}^{K+1}$, where $K$ is the number of
phases in
the service time distribution. Variants of these maps have been
employed in the
literature [see, e.g., \citet{MandelStrong}, \citet{reed2007a},
\citet{DaiTez06}, \citet{daiTez07a} and \citet
{Whitt07c}]. We use these maps not
just in diffusion scaling but also in fluid scaling. Using a map twice, one
for each scaling, allows us to obtain diffusion limits as a simple consequence
of the standard continuous-mapping theorem and the random-time-change theorem.
In the seminal paper, \citet{rei84} proves a \textit{conventional}
heavy traffic
limit theorem for generalized Jackson networks. Our approach resembles the
work of \citet{joh83}, which also uses a multi-dimensional
Skorohod map twice
and provides a significant simplification of Reiman's original proof.

For the $G/\GI/n$ model, \citet{reed2007a} proves a many-server
limit for the
total-customer-count process in the critically-loaded regime; his
assumption on the
service time distribution is completely general, and his limit is not a Markov
process. This work is generalized in \citet
{MandelbaumMomcilovic09} to allow for
customer abandonment. For the overloaded $G/\GI/n$ model, \citet{reed2008}
prove a finite-dimensional-distribution limit for the total-customer-count
process. \citet{JMM04} prove a limit theorem for the $\GI/D/n$ model.
\citet{GamarnikMomcilovic08} study a many-server limit of the steady-state
distribution of the $\GI/\GI/n$ model, where the service times are
lattice-valued with a finite support. When the service time
distribution is
general, measure-valued processes have been used to give a Markovian
description of the system. \citet{kaspi2007} obtain a
measure-valued fluid
limit for the $G/\GI/n$ model. \citet{KangRamanan08} obtain a
measure-valued
fluid limit for the $G/\GI/n+G$ model with customer abandonment, and
\citet{JihengZhang09} obtains a similar measure-valued fluid limit
independently. Their work partially justifies the fluid model in
\citet{Whitt06}.

The remainder of the paper is organized as follows. In Section \ref{sec:GGn},
we introduce the $G/\GI/n+\GI$ model, in an asymptotic framework, and phase-type
distributions. The main results, Theorems \ref{thm:critical}%
--\ref{thm:overloaded}, are stated in Section \ref{sec:queue}; a
roadmap for the
proofs is introduced in Section \ref{sec:roadmap}. In Section \ref
{sec:rep} we
introduce a perturbed system that is equivalent to a $G/\PH/n+M$ queue and
derive the dynamical equations that the perturbed system must obey. The proofs
for the diffusion limits of $G/\PH/n+M$ queues, in both the
critically-loaded and the
overloaded regimes, are given in Section \ref
{sec:exponentialPatience}. Section
\ref{sec:criGI} is dedicated to the proof for the diffusion limit of
$G/\PH/n+\GI$ queues in the critically-loaded regime. In Appendix \ref
{sec:Map}, we
introduce a continuous map and establish various properties for the
map. The
state-space-collapse lemma is proved in Appendix \ref
{sec:proof-theorems}, in
which Theorem \ref{thm:virtualwaiting} and Lemma \ref{lem:aq} are
also proved.

\subsection*{Notation}

All random variables and processes are defined on a common probability space
$(\Omega,\mathcal{F},\mathbb{P})$ unless otherwise specified. The symbols
$\mathbb{Z}$, $\mathbb{Z}_{+}$, $\mathbb{N}$, $\mathbb{R}$ and
$\mathbb
{R}%
_{+}$ are used to denote the sets of integers, nonnegative integers, positive
integers, real numbers and nonnegative real numbers, respectively. For
$d\in\mathbb{N}$, $\mathbb{R}^{d}$ denotes the $d$-dimensional Euclidean
space; thus, $\mathbb{R}=\mathbb{R}^{1}$. The space of functions
$f\dvtx\mathbb{R}_{+}\rightarrow\mathbb{R}^{d}$ that are
right-continuous on
$ [ 0,\infty) $ and have left
limits in $ ( 0,\infty)$ is
denoted by $\mathbb{D}(\mathbb{R}_{+}, \mathbb{R}^{d})$ or simply
$\mathbb{D}^{d}$;
similarly, with $T>0$, the space of functions $f\dvtx [ 0, T ] \rightarrow
\mathbb{R}^{d}$
that are right-continuous on $ [ 0,T ) $ and have left
limits in $ ( 0,T ]$ is denoted by $\mathbb{D}([0,T], \mathbb{R}^{d})$.
For $f\in\mathbb{D}^d$, $f(t-)$ denotes its
left limit at $t>0$. For a sequence of random elements $\{X^n, n\in
\mathbb{N}\}$
taking values in a metric space, we write $X^n\Rightarrow X$ to denote
the convergence of $X^n$ to $X$ in distribution. Each stochastic
process whose
sample paths are in $\mathbb{D}^d$ is considered to be a $\mathbb
{D}^d$-valued random
element. The
space $\mathbb{D}^d$ is assumed to be endowed with the Skorohod $J_1$-topology
[see \citet{ethkur86} or \citet{bil99}].
Given $x\in\mathbb{R}$, we set
$x^{+}=\max\{x,0\}$, $x^{-}=\max\{-x,0\}$ and $ \lfloor x \rfloor
=\max\{j\in\mathbb{Z}\dvtx j\leq x\}$. All vectors are envisioned as column
vectors. For a $K$-dimensional vector $u\in\mathbb{R}^{K}$, we use
$u_{k}$ to denote its $k$th
entry and $\operatorname{diag}(u)$ for the $K\times K$ diagonal
matrix with
$k$th diagonal entry $u_{k}$; we put $ \vert u \vert={\max_{1\leq
k\leq K}} \vert u_{k} \vert$. For a matrix $M$, $M^{\prime}$ denotes
its transpose and $M_{jk}$ denotes its $(j,k)$th entry. We reserve $I$
for the
$K\times K$ identity matrix, $e$ for the $K$-dimensional vector of ones and
$e^{k}$ for the $K$-dimensional vector with $k$th entry one and all other
entries zero.

\section{A $G/\PH/n+\GI$ queue}
\label{sec:GGn}

We first introduce a $G/\GI/n+\GI$ queue in Section \ref{sec:GGnsub}. We then
define a $G/\PH/n+\GI$ queue by restricting the service time distribution to
be phase-type. Phase-type distributions are defined in Section
\ref{sec:phasetype}.

\subsection{A $G/\GI/n+\GI$ queue}
\label{sec:GGnsub}

A $G/\GI/n$ queue is a classical stochastic system that has been extensively
studied in the literature [see, e.g., \citet{KieferWolf55},
\citet{bor67a} and
\citet{iglwhi70a}, among others].
In such a system, there are $n$ identical servers. The service times
$\{v_{i}, i\in\mathbb{N}\}$ are a sequence of
i.i.d. random variables, where $v_{i}$ is the service time of the $i$th customer
entering service after time $0$. The service time distribution is
general (the $\GI$ in the $G/\GI/n$ notation), although for the rest of this
paper we restrict it to be a phase-type distribution. The arrival process
$E=\{E(t),t\geq0\}$ is assumed to be general (the first $G$), where $E(t)$
denotes the number of customer arrivals to the system by time $t$. Upon his
arrival to the system, a customer gets into service immediately if
there is an
idle server; otherwise, he waits in a \textit{waiting buffer} that holds a
first-in-first-out (FIFO) queue. The buffer size is assumed to be infinite.
When a server finishes a service, the server removes the leading
customer from
the waiting buffer and starts to serve the customer; when the queue is empty,
the server begins to idle. In our model, each customer has a patience time:
when a customer's waiting time in queue exceeds his patience time, the
customer abandons the system without any service. Retrial is not
modeled in
this paper. We assume that the patience times of customers
who arrive after time $0$, form a sequence of i.i.d. random variables
that have a general
distribution. Thus our model is a $G/\GI/n+\GI$ queue, where
$+\GI$ signifies the general patience time distribution. When the patience
times are i.i.d. following an exponential distribution, the resulting
system is a
$G/\GI/n+M$ queue.

We focus on systems when the arrival rate is high. Specifically, we
consider a
sequence of $G/\GI/n+\GI$ systems indexed by $n$, with $E^{n}$ being the arrival
process of the $n$th system. We assume that for the $n$th system, the arrival
rate $\lambda^{n}\rightarrow\infty$ as the number of servers
$n\rightarrow
\infty$, while the service time and the patience time distributions do not
change with $n$. We further assume that
%
%
\begin{equation} \label{eq:lambdan}\qquad
\lim_{n\rightarrow\infty}\frac{\lambda^{n}}{n}=\lambda>0,\qquad \lim
_{n\rightarrow\infty}\sqrt{n} \biggl( \lambda- \frac{\lambda^{n}}{n} \biggr)
=\beta\mu\qquad\mbox{for some }\beta\in\mathbb{R}
\end{equation}
and
%
%
\begin{equation}\label{eq:Eclt}%
\tilde{E}^{n}\Rightarrow\tilde{E} \qquad\mbox{as }n\rightarrow\infty,
\end{equation}
where $\tilde{E}$ is a Brownian motion and
%
%
\begin{equation} \label{eq:tilEn}%
\tilde{E}^{n}(t)=\frac{1}{\sqrt{n}}\hat{E}^{n}(t),\qquad \hat{E}^{n}%
(t)=E^{n}(t)-\lambda^{n}t \qquad\mbox{for }t\geq0.
\end{equation}
We use $m$ to denote the mean service time; thus $\mu=1/m$ is the mean service
rate. For future purposes, let
%
%
\begin{equation}
\label{eq:rho}
\rho^{n}=\frac{\lambda^{n}}{n\mu} \quad\mbox{and}\quad
\rho=\frac{\lambda}{\mu}.
\end{equation}
Because customer abandonment is allowed, it is not necessary to assume
$\rho^{n}<1$ or $\rho\leq1$. Condition (\ref{eq:lambdan}) implies that
\[
\lim_{n\rightarrow\infty}\sqrt{n}(\rho-\rho^{n})=\beta.
\]
When $\rho=1$, the sequence of systems is critically loaded in the
limit, and
is said to be in the Quality and Efficiency-Driven (QED) regime or the
Halfin--Whitt regime. When $\rho>1$, the sequence of systems is
overloaded, and
is said to be in the Efficiency-Driven (ED) regime. Our focus is both
the QED
and the ED regimes.

\subsection{Phase-type distributions}
\label{sec:phasetype}

In this section, we introduce phase-type distributions.
Such a distribution is assumed to have $K\geq1$ phases. The set of
phases is
assumed to be $\mathcal{K}=\{1,\ldots,K\}$. Each phase-type
distribution has a
set of parameters $p$, $\nu$ and $P$, where $p$ is a $K$-dimensional vector
of nonnegative entries whose sum is equal to one, $\nu$ is a $K$-dimensional
vector of positive entries and $P$ is a $K\times K$ sub-stochastic
matrix. We
assume that the diagonals of $P$ are zero, namely,
%
%
\begin{equation} \label{eq:Piizero}%
P_{ii}=0 \qquad\mbox{for }i=1,\ldots,K,
\end{equation}
and $P$ is transient, namely,
%
%
\begin{equation} \label{eq:transient}%
I-P\qquad\mbox{is invertible}.
\end{equation}
A (continuous) phase-type random variable $v$ is defined as the time until
absorption in an absorbing state of a continuous-time Markov chain.
With $p$,
$\nu$ and $P$, the continuous-time Markov chain can be described as follows.
It has $K+1$ states, $1,\ldots,K,K+1$, with state $K+1$ being
absorbing. The
rate matrix (or generator) of the Markov chain is
\[
G=%
\pmatrix{
F & h\cr
0 & 0
}
,
\]
where $F=\operatorname{diag}(\nu)(P-I)$ is a $K\times K$ matrix and
$h=-Fe$ is a
$K$-dimensional vector.
\begin{definition}
A continuous phase-type random variable $v$ with parameters $p$, $\nu$ and
$P$, denoted as $v\sim\PH(p, \nu, P)$, is defined to
be the
first time until the continuous-time Markov chain with initial distribution
$p$ and rate matrix $G$ reaches state $K+1$.
\end{definition}

Given condition (\ref{eq:Piizero}), the rate matrix $G$ and $(\nu,P)$ are
uniquely determined from each other. It is well known [see, e.g.,
\citet{LatoucheRamaswami}] that
\[
\mathbb{P}[v\leq x]=1-p^{\prime}\exp(Fx)e \qquad\mbox{for }x\geq0.
\]
Because parameters $p$, $\nu$ and $P$ uniquely determine a phase-type
distribution, we free symbols $F$ and $G$ so that they can be reused in the
rest of the paper.

For our purposes, we provide an alternative way to sample a
$\PH(p,\nu,P)$ random variable. We first sample a
sequence of
phases $k_{1},\ldots,k_{L}$ in $\mathcal{K}=\{1,\ldots,K\}$ as
follows. We
sample phase $k_{1}$ following distribution $p$ on $\mathcal{K}$. Assume
$k_{1},\ldots,k_{i}\in\mathcal{K}$ have been sampled. Setting $j=k_{i}$,
sample a phase from $\{1,\ldots,K,K+1\}$ following a distribution that is
determined by the $j$th row of $P$; the probability of getting phase
$K+1$ is
$1-\sum_{\ell=1}^{K}P_{j\ell}\geq0$. The resulting phase is denoted by
$k_{i+1}$. When $k_{i+1}=K+1$, terminate the process and set $L=i$; otherwise,
continue the sampling process. Because the matrix $P$ is assumed to be
transient, one has $L<\infty$ almost surely. Let $\xi_{1},\ldots,\xi
_{L}$ be
independently sampled from exponential distributions with respective
rates $\nu_{k_{1}%
},\ldots,\nu_{k_{L}}$. Then%
%
%
\begin{equation} \label{eq:phasev}%
v=\sum_{i=1}^{L}\xi_{i}.
\end{equation}

\section{Main results}
\label{sec:queue}

In this section, we present two sets of results. The first set,
presented in
Section \ref{sec:GPhn+GI}, is for critically loaded $G/\PH/n+\GI$
queues. The
second set, presented in Section \ref{sec:diff-limit-overl}, is for overloaded
$G/\PH/n+M$ queues. A roadmap for proving these results is given in
Section \ref{sec:roadmap}.

We consider a sequence of $G/\PH/n+\GI$ queues, indexed by the number of servers
$n$, which satisfies condition (\ref{eq:lambdan}). We assume that the service
times follow a phase-type distribution $\PH
(p,\nu,P)$. Each customer's service time can be decomposed into a
number of
phases as in (\ref{eq:phasev}). When a customer is in service, it must
be in
one of the $K$ phases of service. Let $Z_{k}^{n}(t)$ denote the number of
customers \textit{in phase $k$ service} in the $n$th system at time $t$; service
times in phase $k$ are exponentially distributed with rate $\nu_{k}$.
We use
$Z^{n}(t)$ to denote the corresponding $K$-dimensional vector. We call
$Z^{n}=\{Z^{n}(t),t\geq0\}$ the \textit{server-allocation process}. Let
$N^{n}(t)$ denote the number of customers in the $n$th system at time $t$,
either in queue or in service. Setting
%
%
\begin{equation} \label{eq:X}%
X^{n}(t)=N^{n}(t)-n \qquad\mbox{for }t\geq0,
\end{equation}
we call $X^{n}=\{X^{n}(t),t\geq0\}$ the \textit{total-customer-count
process} in
the $n$th system. One can check that $(X^{n}(t))^{+}$ is the number of
customers waiting in queue at time~$t$, and $(X^{n}(t))^{-}$ is the
number of
idle servers at time $t$. Clearly,
%
%
\begin{equation} \label{eq:XZ}%
e^{\prime}Z^{n}(t)=n-(X^{n}(t))^{-} \qquad\mbox{for }t\geq0.
\end{equation}
The processes $X^{n}$ and $Z^{n}$ describe the ``state'' of the system as time evolves. Hereafter, they are
called the \textit{state processes} for the $n$th system.

The customers in service are distributed among the $K$ phases following a
distribution $\gamma$, given by
%
%
\begin{eqnarray}
\label{eq:gamma}
\gamma&=&\mu R^{-1}p,\\
\label{eq:R}%
R&=&(I-P^{\prime})\operatorname{diag}(\nu).
\end{eqnarray}
One can check that $\sum_{k=1}^{K}\gamma_{k}=1$, and one interprets
$\gamma_{k}$ to be the fraction of phase $k$ load on the $n$ servers.

The preceding paragraph suggests the following centering for the
server-allocation process:
\[
\hat{Z}^{n}(t)=Z^{n}(t)-n\gamma\qquad\mbox{for }t\geq0.
\]
Define the corresponding diffusion-scaled process
\[
\tilde{Z}^{n}(t)=\frac{1}{\sqrt{n}}\hat{Z}^{n}(t) \qquad\mbox{for }t\geq0.
\]

\subsection{Diffusion limits for critically loaded $G/\PH/n+\GI$ queues}
\label{sec:GPhn+GI}

Throughout Section \ref{sec:GPhn+GI}, we assume that
$\rho=1$ and that the patience times of customers who arrive after
time $0$
are i.i.d. having distribution function $F$, which satisfies
%
%
\begin{equation} \label{1}%
F(0)=0 \quad\mbox{and}\quad \alpha=\lim_{x\downarrow0}x^{-1}F ( x )
<\infty.
\end{equation}
Note that for exponentially distributed patience times, $\alpha$ turns
out to
be the rate of the exponential distribution.

Since $\rho$ is assumed to be $1$, we define the diffusion-scaled
total-customer-count process $\tilde{X}^{n}$ by%
%
%
\begin{equation}\label{eq:diffusionScaleCri}%
\tilde{X}^{n}(t)=\frac{1}{\sqrt{n}}X^{n}(t)\qquad \mbox{for }t\geq0.
\end{equation}
[When $\rho>1$, the definition of $\tilde{X}^{n}(t)$ will be
modified, which
is given in (\ref{eq:XhatOL}).] We assume that
%
%
\begin{equation}\label{eq:Zinitial}%
(\tilde{X}^{n}(0),\tilde{Z}^{n}(0))\Rightarrow(\tilde{X}(0),\tilde
{Z}%
(0)) \qquad\mbox{as }n\rightarrow\infty
\end{equation}
for a pair of random variables $(\tilde{X}(0),\tilde{Z}(0))$.

The random variables $\tilde{X}(0)$ and $\tilde{Z}(0)$ are assumed to be
defined on some probability space $(\tilde{\Omega},\tilde{\mathcal
{F}}%
,\tilde{\mathbb{P}})$, which is rich enough so that stochastic processes
$\tilde{E},\tilde{\Phi}^{0},\ldots,\tilde{\Phi}^{K}$ and $\tilde
{S}$ defined
on this space are independent of $(\tilde{X}(0),\tilde{Z}(0))$. Here,
$\tilde{E}$ is a one-dimensional driftless Brownian motion, and
$\tilde
{\Phi
}^{0},\ldots,\tilde{\Phi}^{K}$ and $\tilde{S}$ are $K$-dimensional driftless
Brownian motions. These Brownian motions, possibly degenerate, are mutually
independent and start from $0$; the variance of $\tilde{E}$ is
$\lambda
c_{a}^{2}$ for some constant $c_{a}^{2}\geq0$, and the covariance
matrices of
$\tilde{\Phi}^{0},\ldots,\tilde{\Phi}^{K}$ and $\tilde{S}$ are
$H^{0}%
,\ldots,H^{K}$ and $\operatorname{diag}(\nu)$, respectively, where
for $k=0,\ldots
,K$, the
$K\times K$ matrix $H^{k}$ is given by
%
%
\begin{equation}\label{eq:Hk}%
H_{ij}^{k}=%
\cases{
p_{i}^{k}(1-p_{j}^{k}), &\quad if $i=j$,\cr
-p_{i}^{k}p_{j}^{k}, &\quad otherwise,}
\end{equation}
with $p^{0}=p$ and $p^{k}$ being the $k$th column of $P^{\prime}$.

To state the main theorems of this paper, let
%
%
\begin{eqnarray}
\label{eq:tU}
\tilde{U}(t)&=&\tilde{X}(0)+\tilde{E}(t)-\mu\beta t+e^{\prime
}\tilde
{M}(t),\\
\label{eq:tV}%
\tilde{V}(t)&=&(I-pe^{\prime})\tilde{Z}(0)+\tilde{\Phi}^{0}(\mu
t)+(I-pe^{\prime})\tilde{M}(t),
\end{eqnarray}
where
%
%
\begin{equation} \label{eq:tM}%
\tilde{M}(t)=\sum_{k=1}^{K}\tilde{\Phi}^{k}(\nu_{k}\gamma
_{k}t)-(I-P^{\prime
})\tilde{S}(\gamma t)
\end{equation}
for $t\geq0$. The process $(\tilde{U},\tilde{V})$ is a
$(K+1)$-dimensional Brownian motion; it is degenerate because $e'
\tilde V(t)$ = 0 for $t \geq0$. [When $\rho>1$, the definition of
$\tilde{U}$ will be modified in (\ref{eq:tUOvr}).] Before we state
the first theorem of this paper, we present the following lemma,
which is a corollary of Lemma
\ref{lem:maps} in Appendix \ref{sec:Map}.

\begin{lemma}
\label{lem:Phi}Let $p$ be a $K$-dimensional vector that is the
distribution of
initial phases of the phase-type service times, $R$ be the $K\times K$ matrix
defined by (\ref{eq:R}) and $\alpha\geq0$ be defined by (\ref{1}).
\textup{(a)} For
each $(u,v)\in\mathbb{D}^{K+1}$ with $u(t)\in\mathbb{R}$ and
$v(t)\in
\mathbb{R}^{K}$ for $t\geq0$,
there exists a unique
$(x,z)\in\mathbb{D}^{K+1}$ with $x(t)\in\mathbb{R}$ and $z(t)\in
\mathbb
{R}^{K}$ for $t\geq0$, such that
%
%
\begin{eqnarray}
\label{eq:mapPhi1}
x(t)&=&u(t)-\alpha\int_{0}^{t}(x(s))^{+} \,ds-e^{\prime}R\int_{0}%
^{t}z(s)\, ds,\\
\label{eq:mapPhi2}%
z(t)&=&v(t)-p(x(t))^{-}-(I-pe^{\prime})R\int_{0}^{t}z(s)\, ds
\end{eqnarray}
for $t\geq0$. \textup{(b)} For each $(u,v)\in\mathbb{D}^{K+1}$, define
$\Phi(u,v)=(x,z)\in\mathbb{D}^{K+1}$, where $(x,z)$ satisfies
(\ref{eq:mapPhi1}) and (\ref{eq:mapPhi2}). The map $\Phi$ is well
defined and
is continuous when both the domain and the range $\mathbb{D}^{K+1}$
are endowed with the Skorohod $J_{1}$-topology. \textup{(c)}
The map $\Phi$ is Lipschitz continuous in the sense that for any
$T>0$, there
exists a constant $C_{T}^{1}>0$ such that
%
%
\begin{eqnarray} \label{eq:LipschitzPhi}%
{\sup_{0\leq t\leq T}} \vert\Phi(y)(t)-\Phi(\tilde{y})(t) \vert\leq
{C_{T}^{1}\sup_{0\leq t\leq T}} \vert y(t)-\tilde{y}(t) \vert
\nonumber\\[-8pt]\\[-8pt]
\eqntext{\mbox{for any }y,\tilde{y}\in\mathbb{D}^{K+1}.}
\end{eqnarray}
\textup{(d)} The map $\Phi$ is positively homogeneous in the sense that
%
%
\begin{equation} \label{eq:PhiPhomogenous}%
\Phi(ay)=a\Phi(y) \qquad\mbox{for each }a>0\mbox{ and each }y\in\mathbb
{D}%
^{K+1}.
\end{equation}
\end{lemma}

Let $A_{0}^{n}$ be the number of customers who are waiting in queue at time
$0$ but will eventually abandon the system, and
\[
\tilde{A}_{0}^{n}=\frac{1}{\sqrt{n}}A_{0}^{n}.
\]
To state Theorem \ref{thm:critical}, we assume that
%
%
\begin{equation}\label{eq:initial}%
\tilde{A}_{0}^{n}\Rightarrow0 \qquad\mbox{as }n\rightarrow\infty.
\end{equation}
Clearly, condition (\ref{eq:initial}) is satisfied if no customers are
waiting in queue at time~$0$. The validity of this initial condition
will be further discussed at the end of Section~\ref{sec:GPhn+GI}.

\begin{theorem}
\label{thm:critical} Consider a sequence of $G/\PH/n+\GI$ queues satisfying
(\ref{eq:lambdan}) and (\ref{eq:Eclt}). Assume that $\rho=1$ and that
(\ref{1}), (\ref{eq:Zinitial}) and (\ref{eq:initial}) hold. Then
\[
(\tilde{X}^{n},\tilde{Z}^{n})\Rightarrow(\tilde{X},\tilde{Z})
\qquad\mbox{as
}n\rightarrow\infty,
\]
where
%
%
\begin{equation} \label{eq:criX}%
(\tilde{X},\tilde{Z})=\Phi(\tilde{U},\tilde{V}).
\end{equation}
\end{theorem}

Suppose that each customer, including those initial customers
who are waiting in queue at time $0$, samples his first service phase
that he is yet to enter following
distribution $p$ at his arrival time to the system. One can stratify the
customers in the waiting buffer according to their first service
phases. For
$k=1,\ldots,K$, we use $Q_{k}^{n}(t)$ to denote the number of waiting
customers at time $t$ whose initial service phase will be phase $k$
[$Q_{k}%
^{n}(t)=0$ for $t\geq0$ if phase $k$ is not a first service phase for any
customer], and we use $Y_{k}^{n}(t)$ to denote the number of phase $k$
customers in the system at time $t$, either waiting or in service. Let
$Q^{n}(t)$ and $Y^{n}(t)$ denote the corresponding $K$-dimensional vectors.
Set%
%
%
\begin{eqnarray} \label{eq:VectorQ}%
\tilde{Q}_{k}^{n}(t)&=&\frac{1}{\sqrt{n}}Q_{k}^{n}(t),\qquad \tilde{Y}_{k}%
^{n}(t)=\frac{1}{\sqrt{n}}\hat{Y}_{k}^{n}(t),\nonumber\\[-8pt]\\[-8pt]
\hat{Y}_{k}^{n}(t)&=&Y_{k}^{n}(t)-n\gamma_{k} \qquad\mbox{for
}t\geq0.\nonumber
\end{eqnarray}
Clearly,
%
%
\begin{equation}\label{eq:XY}%
\tilde{Y}_{k}^{n}(t)=\tilde{Q}_{k}^{n}(t)+\tilde{Z}_{k}^{n}(t) \quad\mbox{and}\quad
\tilde{X}^{n}(t)=e^{\prime}\tilde{Y}^{n}(t) \qquad\mbox{for }t\geq0.
\end{equation}

The following lemma says that for critically loaded systems, the waiting
customers are distributed among the $K$ phases following distribution
$p$. It
is known as the state-space-collapse (SSC) result.
\begin{lemma}
\label{thm:ssc} Under the conditions of Theorem \ref{thm:critical},
for any
$T>0$,
%
%
\begin{equation}\label{eq:SSC}%
{\frac{1}{\sqrt{n}}\sup_{0\leq t\leq T}} \vert Q^{n}(t)-p(X^{n}%
(t))^{+} \vert\Rightarrow0 \qquad\mbox{as }n\rightarrow\infty.
\end{equation}
\end{lemma}

The following theorem is a corollary to Theorem \ref{thm:critical} and
Lemma \ref{thm:ssc}. When there is no customer abandonment and the arrival
process is renewal, it is identical to Theorem 2.3 of \citet
{PuhalskiiReiman00}.
\begin{theorem}
\label{cor:1} Under the conditions of Theorem \ref{thm:critical},
\[
(\tilde{X}^{n},\tilde{Y}^{n},\tilde{Z}^{n})\Rightarrow(\tilde
{X},\tilde
{Y},\tilde{Z}) \qquad\mbox{as }n\rightarrow\infty,
\]
where $(\tilde{X},\tilde{Z})$ is defined in (\ref{eq:criX}) and
%
%
\begin{equation}\label{eq:tF}%
\tilde{Y}(t)=p(\tilde{X}(t))^{+}+\tilde{Z}(t).
\end{equation}
The process $\tilde{Y}$ satisfies
\begin{eqnarray*}
\tilde{Y}(t) & = & \tilde{Y}(0)-\beta\mu pt+\tilde{\Phi}^{0}(\mu
t)+p\tilde
{E}(t)+\tilde{M}(t)\\
&&{} - R\int_{0}^{t}\tilde{Y}(s) \,ds+(R-\alpha I)p\int_{0}^{t}(e^{\prime
}\tilde{Y}(s))^{+} \,ds \qquad\mbox{for $t\geq0$}.
\end{eqnarray*}
\end{theorem}

The process $\tilde{Y}$ in Theorem \ref{cor:1} is a diffusion process (see
\citeauthor{RogersWilliams2} [(\citeyear{RogersWilliams2}), page 110]
or \citeauthor{karlintaylor1} [(\citeyear{karlintaylor1}), page
159] for a definition of diffusion processes). Therefore,
$\tilde{Y}$ is a continuous Markov process. The map $\Phi$ in (\ref{eq:criX})
defines $(\tilde{X},\tilde{Z})$ as a $(K+1)$-dimensional continuous process,
which is degenerate because it lives on a $K$-dimensional manifold.
From the
$K$-dimensional process $\tilde{Y}$, one can recover the $(K+1)$-dimensional
process $(\tilde{X},\tilde{Z})$ via
%
%
\begin{equation}\label{eq:XZfromY}%
\tilde{X}(t)=e^{\prime}\tilde{Y}(t) \quad\mbox{and}\quad \tilde{Z}%
(t)=\tilde{Y}(t)-p(\tilde{X}(t))^{+} \qquad\mbox{for $t\geq0$}.
\end{equation}
Therefore, $(\tilde{X},\tilde{Z})$ is also a continuous Markov process.
However, the process $(\tilde{X},\tilde{Z})$ is not a diffusion process
by the
common definition because the function $x^{+}$ in (\ref{eq:XZfromY})
is not
twice differentiable in $x$ at $0$. \citeauthor{Whitt05a} [(\citeyear
{Whitt05a}), Remark 2.2] makes a
similar observation that his limit process is not a diffusion process,
but a
simple transformation of his limit process is a diffusion process.

Our next theorem is concerned with the \textit{virtual waiting time process}
$W^{n}=\{W^{n}(t), t\ge0\}$. Here, $W^{n}(t)$ is the potential waiting
time of
a hypothetical, infinitely patient customer who arrives at the queue at time
$t$. When there is no customer abandonment and the arrival process is renewal,
the theorem is implied by Corollary 2.3 and Remark 2.6 of
\citet{PuhalskiiReiman00}.

\begin{theorem}
\label{thm:virtualwaiting} Under the conditions of Theorem \ref
{thm:critical},
%
%
\begin{equation} \label{eq:Vw}%
\sqrt{n}W^{n}\Rightarrow\frac{(\tilde{X})^{+}}{\mu} \qquad\mbox{as
}n\rightarrow\infty.
\end{equation}
\end{theorem}

We end this section by the following lemma, which gives a
justification for imposing initial condition (\ref{eq:initial}) in
Theorems \ref{thm:critical}--\ref{thm:virtualwaiting}.
Let $A_{Q}^{n}(t)$ be the number of customers in the $n$th system who
are waiting in queue at time $t$, but will eventually abandon the system.
Clearly,
\[
A_{0}^{n}=A_{Q}^{n}(0).
\]
Its diffusion-scaled version is given by%
%
%
\begin{equation}\label{eq:AQTilde}%
\tilde{A}_{Q}^{n}(t)=\frac{1}{\sqrt{n}}A_{Q}^{n}(t)\qquad\mbox{for }%
t\geq0.
\end{equation}
Regarding the process $\tilde{A}_{Q}^{n}=\{\tilde{A}_{Q}^{n}(t),t\geq
0\}
$, we
have the following result.
\begin{lemma}
\label{lem:aq}Under the conditions of Theorem \ref{thm:critical},%
%
%
\begin{equation}\label{eq:AQ}%
\tilde{A}_{Q}^{n}\Rightarrow0 \qquad\mbox{as }n\rightarrow\infty.
\end{equation}
\end{lemma}

The proof of Lemma \ref{lem:aq} is presented in Appendix
\ref{sec:proof-theorems}. Assume that the queue is initially empty.
Then condition (\ref{eq:initial}) is satisfied at time $t=0$.
Under an additional assumption (\ref{eq:Zinitial}), Theorem
\ref{thm:critical} and Lemma \ref{lem:aq} imply that for any $t>0$,
\[
\tilde{A}_{Q}^{n}(t)\Rightarrow0 \quad\mbox{and}\quad
(\tilde X^n(t), \tilde Z^n(t)) \Rightarrow(\tilde X(t), \tilde Z(t))
\qquad\mbox{as } n \rightarrow\infty.
\]
Thus, if we start to observe the system at any fixed time $t>0$, initial
conditions (\ref{eq:Zinitial}) and (\ref{eq:initial}) are
indeed satisfied at time $t$. Condition (\ref{eq:initial}) is used to
prove asymptotic relationship (\ref{eq:DaiHe09}) in the critically-loaded
regime; this relationship between the abandonment-count process and
the queue-length process is the key to extending the diffusion limits
for $G/\PH/n+M$ queues to the $G/\PH/n+\GI$ model with a general patience
time distribution. Condition (\ref{eq:initial}) is necessary for the
asymptotic relationship to hold. In \citet
{MandelbaumMomcilovic09}, an
initial assumption similar to (\ref{eq:initial}) is made for the
$G/\GI/n+\GI$ model in the critically-loaded regime.

\subsection{A diffusion limit for overloaded $G/\PH/n+M$ queues}
\label{sec:diff-limit-overl}

Our next result is for overloaded $G/\PH/n+M$ systems, where
the patience times of all customers, including those waiting in queue
at time $0$,
are assumed to be i.i.d. following an exponential distribution. We use
$\alpha$
to denote the rate of the exponential patience time distribution. Note that
this definition of $\alpha$ is consistent with the definition in~(\ref{1}).
Assume that $\rho>1$. Intuitively, when $n$ is large, all $n$ servers are
$100\%$ busy, and there should be $nq$ customers on average waiting in the
buffer, where
%
%
\begin{equation} \label{eq:q}%
q=\frac{\lambda-\mu}{\alpha}.
\end{equation}
An intuitive explanation is as follows: $\lambda-\mu$ is the number of
customers per unit of time that the system must ``delete''
in order for the system to reach an equilibrium.
While in equilibrium, each waiting customer abandons the system at rate
$\alpha$, and collectively all $q$ waiting customers abandon the
system at
rate $q\alpha$ customers per unit of time. Thus, one should have
$q\alpha=\lambda-\mu$, which leads to (\ref{eq:q}). Readers are
referred to
\citet{Whitt04a} for further discussion on the derivation of
(\ref{eq:q}).

Now we modify the definition of $\tilde{X}^{n}$ in (\ref
{eq:diffusionScaleCri}%
) and $\tilde{U}$ in (\ref{eq:tU}). Let
%
%
\begin{eqnarray}
\label{eq:XhatOL}
\tilde{X}^{n}(t)&=&\frac{1}{\sqrt{n}}\hat{X}^{n}(t),\qquad \hat{X}^{n}%
(t)=X^{n}(t)-nq,\\
\label{eq:tUOvr}%
\tilde{U}(t)&=&\tilde{X}(0)+\tilde{E}(t)-\mu\beta t+e^{\prime
}\tilde
{M}(t)-\tilde{G}(qt)
\end{eqnarray}
for $t\geq0$. In (\ref{eq:tUOvr}), the process $\tilde{G}=\{\tilde
{G}%
(t),t\geq0\}$ is a one-dimensional driftless Brownian motion starting from
$0$, which has variance $\alpha$ and is independent of $\tilde
{E},\tilde
{\Phi
}^{0},\ldots,\tilde{\Phi}^{K}$ and $\tilde{S}$ [recall that $\tilde
{E}%
,\tilde{\Phi}^{0},\ldots,\tilde{\Phi}^{K}$ and $\tilde{S}$ are Brownian
motions defined in Section \ref{sec:GPhn+GI}, and $\tilde{M}$ is
given by
(\ref{eq:tM})]. When $\rho=1$, one has $q=0$. Thus, definitions in
(\ref{eq:XhatOL}) and (\ref{eq:tUOvr}) are consistent with
(\ref{eq:diffusionScaleCri}) and (\ref{eq:tU}). Assume that
%
%
\begin{equation}\label{eq:OLinitial}%
(\tilde{X}^{n}(0),\tilde{Z}^{n}(0))\Rightarrow(\tilde{X}(0),\tilde
{Z}%
(0)) \qquad\mbox{as }n\rightarrow\infty
\end{equation}
for a pair of random variables $(\tilde{X}(0),\tilde{Z}(0))$.

Before presenting Theorem \ref{thm:overloaded}, we introduce the
next lemma, which is also a corollary of Lemma \ref{lem:maps}.
\begin{lemma}
\label{lem:Psi}Let $p$ be a $K$-dimensional vector that is the
distribution of
initial phases of the phase-type service times, $R$ be the $K\times K$ matrix
defined by (\ref{eq:R}) and $\alpha\geq0$ be defined by (\ref{1}).
\textup{(a)} For
each $(u,v)\in\mathbb{D}^{K+1}$ with $u(t)\in\mathbb{R}$ and
$v(t)\in
\mathbb{R}^{K}$ for $t\geq0$,
there exists a unique
$(x,z)\in\mathbb{D}^{K+1}$ with $x(t)\in\mathbb{R}$ and $z(t)\in
\mathbb
{R}^{K}$ for $t\geq0$, such that
%
%
\begin{eqnarray}
\label{eq:mapPsi1}
x(t)&=&u(t)-\alpha\int_{0}^{t}x(s) \,ds-e^{\prime}R\int_{0}^{t}%
z(s) \,ds,\\
\label{eq:mapPsi2}%
z(t)&=&v(t)-(I-pe^{\prime})R\int_{0}^{t}z(s) \,ds
\end{eqnarray}
for $t\geq0$. \textup{(b)} For each $(u,v)\in\mathbb{D}^{K+1}$, define
$\Psi(u,v)=(x,z)\in\mathbb{D}^{K+1}$, where $(x,z)$ satisfies
(\ref{eq:mapPsi1}) and (\ref{eq:mapPsi2}). The map $\Psi$ is well
defined and
is continuous when both the domain and the range $\mathbb{D}^{K+1}$ are
endowed with the Skorohod $J_{1}$-topology. \textup{(c)}
The map $\Psi$ is Lipschitz continuous in the sense that for any
$T>0$, there
exists a constant $C_{T}^{2}>0$ such that
%
%
\begin{eqnarray} \label{eq:LipschitzPsi}
{\sup_{0\leq t\leq T}} \vert\Psi(y)(t)-\Psi(\tilde{y})(t) \vert\leq
{C_{T}^{2}\sup_{0\leq t\leq T}} \vert y(s)-\tilde{y}(s)
\vert\nonumber\\[-8pt]\\[-8pt]
\eqntext{\mbox{for any }y,\tilde{y}\in\mathbb{D}^{K+1}.}
\end{eqnarray}
\textup{(d)} The map $\Psi$ is positively homogeneous in the sense that
%
%
\begin{equation} \label{eq:PsiPhomogenous}%
\Psi(ay)=a\Psi(y) \qquad\mbox{for each }a>0\mbox{ and each }y\in\mathbb
{D}%
^{K+1}.
\end{equation}
\end{lemma}
\begin{theorem}
\label{thm:overloaded}Consider a sequence of $G/\PH/n+M$ queues satisfying
(\ref{eq:lambdan}) and (\ref{eq:Eclt}). Assume that $\rho>1$ and that
(\ref{eq:OLinitial}) holds. Then%
\[
(\tilde{X}^{n},\tilde{Z}^{n})\Rightarrow(\tilde{X},\tilde{Z})
\qquad\mbox{as }n\rightarrow\infty,
\]
where
%
%
\begin{equation} \label{eq:overlNZ}%
(\tilde{X},\tilde{Z})=\Psi(\tilde{U},\tilde{V}).
\end{equation}
\end{theorem}

Equation (\ref{eq:overlNZ}) defines $(\tilde{X},\tilde{Z})$ as a
$(K+1)$-dimensional diffusion process, which is also degenerate and
lives on a
$K$-dimensional manifold.

\subsection{A roadmap for proofs}
\label{sec:roadmap}

Theorems \ref{thm:critical} and \ref{thm:overloaded} are the main
results of this paper. Theorem \ref{thm:overloaded} and a restricted
version of Theorem \ref{thm:critical} are proved in
Section \ref{sec:exponentialPatience}; the restriction is to assume
that the patience times are exponentially distributed. These proofs
use standard FCLTs and then apply the continuous-mapping theorem and
the random-time-change theorem. To construct appropriate continuous
maps, we introduce a perturbed system in Section
\ref{sec:perturbed}, which is equal in distribution to the original
system when the patience time distribution is exponential. Using the
perturbed system, we are able to construct a set of system equations
in Section \ref{sec:dynamics}, which is critical to define the
continuous maps.

Section \ref{sec:criGI} is devoted to proving the general version of Theorem
\ref{thm:critical}. When the patience time distribution is general and the
systems are critically loaded, we first modify the preceding system equations
slightly by replacing the cumulative number of customer abandonments by an
integral of the queue-length process. We then apply an asymptotic relationship
in \citet{DaiHe09} to establish a result that the error from the
replacement is
negligible under a stochastic boundedness assumption of the diffusion-scaled
queue-length process; the latter assumption holds by a comparison
result in
\citet{DaiHe09} and the restricted version of Theorem \ref
{thm:critical} proved
in Section \ref{sec:exponentialPatience}.

Theorem \ref{cor:1} is a corollary to Theorem \ref{thm:critical} and Lemma
\ref{thm:ssc}; the latter is proved in Appendix \ref{sec:proof-theorems}.
Theorem \ref{thm:virtualwaiting} is also proved in Appendix
\ref{sec:proof-theorems}.

\section{System representation for a $G/\PH/n+M$ queue}
\label{sec:rep}

In this section, we first describe a perturbed system
of the
$G/\PH/n+M$ model, and show that this perturbed system is equivalent to the
$G/\PH/n+M$ queue. We then develop the dynamical equations that the perturbed
system must satisfy.

\subsection{A perturbed system}
\label{sec:perturbed}

Now we describe a perturbation of the $G/\PH/n+M$ model. In the perturbed
system, each phase has a \textit{service queue} for the customers ``in
service.'' Only the leading customer in the service queue is actually in
service; all others are waiting in the service queue, ordered according
to the
FIFO discipline. We use $Z^{*}_{k}(t)$ to denote the number of
customers in
phase $k$ service queue at time $t$. (Star-version quantities are associated
with the perturbed system; the corresponding quantities in the original system
are denoted by the same symbols without the star.) Each customer in the
service queue is attached to exactly one server. Thus, there are exactly
$Z^{*}_{k}(t)$ servers in phase $k$ at time $t$. All these $Z^{*}_{k}(t)$
servers simultaneously work on the leading customer in the service
queue. The
service effort received by the leading customer is additive,
proportional to
the number of servers working on the customer. Each customer has a service
requirement for each phase that he visits, with phase $k$ service requirement
being exponentially distributed with mean $1/\nu_{k}$.

When the total service effort spent on a customer reaches his service
requirement in a phase, the service in the phase is completed. When a customer
completes a phase $k$ service, he immediately moves to the next phase
following a sampling procedure to be specified below, taking with him the
associated server. If the service queue in the new phase is not empty
at their
arrival, the server joins the service immediately, collaborating with other
servers who are already in service to work on the leading customer in the
service queue. The newly arrived customer joins the end of the service queue.
If the new service queue is empty, the server works on her customer who
is the
only one in the new phase of service.

When a customer finishes a phase $k$ service, it uses the $k$th row of
$P$ to
sample the next phase of service to join among phases $\{1, \ldots, K,
K+1\}$;
the probability of selecting phase $K+1$ is $1-\sum_{\ell=1}^{K} P_{k
\ell}$.
If $\ell\in\{1, \ldots, K\}$ is selected, both the server and the customer
move next to phase $\ell$. If $K+1$ is selected, the customer exits
the system
and the associated server is released. The released server checks the FIFO
real queue to select the next customer to work on if the real queue is not
empty. The selected customer is attached to the server until the customer
exits the system. If the real queue is empty, the server becomes idle.

At a customer's arrival time to the system, if there is an idle server, the
customer grabs a free server and is attached to the server until the customer
exits the system. Together, they move into the customer's first phase of
service, which is selected according to distribution $p$. The service and
waiting mechanism is identical to the one described in the previous paragraph.
If all servers are busy at the customer's arrival time, the customer
joins the
end of the FIFO real queue. Only the leading customer in the FIFO real queue
can abandon the system; other waiting customers are infinitely patient until
they become a leading customer. The patience time of the leading
customer is
exponentially distributed with mean $1/\alpha$. The customer abandons the
system without service if his patience clock exceeds the patience time. The
patience clock starts from $0$ when the customer becomes a leading customer
and increases at rate $k$ when the queue length is $k$.

For each $n$ fixed, now we show that when the arrival process $E^{n}$
is a
renewal process, the perturbed system and the original system are equivalent
in a precise mathematical sense. For that, recall that the waiting
buffer in
the perturbed system maintains a FIFO queue for waiting customers. Let
\[
\mathcal{Q}^{\ast}(t)=\bigl(i_{1},\ldots,i_{L^{\ast}(t)}\bigr),
\]
where $L^{\ast}(t)$ is the total number of customers waiting in queue
at time
$t$, and $i_{\ell}$ is the \textit{first service phase} that the $\ell$th
customer is yet to enter.

Let $\xi(t)$ be the remaining interarrival time at time $t$. ($\xi$
has no
star because the arrival processes in the perturbed system and in the original
system are identical.) It follows that, $\{(\xi(t),\mathcal{Q}^{\ast
}(t),Z^{\ast}(t)),t\geq0\}$ is a continuous-time Markov process
living in
state space $\mathbb{R}_{+}\times\mathcal{K}^{\infty}\times\mathbb
{Z}_{+}^{K}%
$, where $\mathcal{K}^{\infty}$ is the space of finite sequences
taking values
in $\mathcal{K}=\{1,\ldots,K\}$.

Let $\{(\xi(t),\mathcal{Q}(t),Z(t)),t\geq0\}$ be the corresponding
process of
the original system. The process is also a continuous-time Markov
process. At
any time $t$, the phase $k$ service rate is $Z_{k}^{\ast}(t)\nu_{k}$
in the
perturbed system and $Z_{k}(t)\nu_{k}$ in the original system, while the
abandonment rate is $L^{\ast}(t)\alpha$ in the perturbed system and
$L(t)\alpha$ in the original system. One can check that the two Markov
processes
\[
\{(\xi(t),\mathcal{Q}(t),Z(t)),t\geq0\}\quad \mbox{and}\quad \{(\xi
(t),\mathcal{Q}^{\ast}(t),Z^{\ast}(t)),t\geq0\}
\]
have the same generator. Thus, when they have the same initial distribution,
they have the same distribution for the entire processes. In the
following, we
always choose the initial condition of the perturbed system to be
identical to
that of the original system.

Even if the arrival process is not a renewal process, the perturbed
system can
still have the same distribution as the original system. The rest of
the paper
does not require the arrival process to be renewal. Rather, we assume that
each arrival process satisfies the requirement that the perturbed
system has
the same distribution as the original one. See \citet{Tezcan06}
for a more
general treatment of perturbed systems.

\subsection{System equations}
\label{sec:dynamics}

From now on, we focus on the perturbed system of the
$G/\PH/n+M$ queue and drop the stars attached to its quantities. We assume
that the patience times of all customers, including those who are
waiting in queue at time $0$, are exponentially distributed with rate
$\alpha$.
In this section, we describe the dynamical equations that the system
must obey. For
$k=1,\ldots,K$, let $\phi^{k}=\{\phi^{k}(j),j\in\mathbb{N}\}$ be a
sequence of
i.i.d. ``Bernoulli random vectors.''
For each
$j$, the $K$-dimensional random vector $\phi^{k}(j)$ takes vector
$e^{\ell}$
with probability $P_{k\ell}$ and takes the $K$-dimensional zero vector with
probability $1-\sum_{\ell=1}^{K}P_{k\ell}$. Similarly, let $\phi
^{0}=\{
\phi^{0}(j),j\in\mathbb{N}\}$ be
a sequence of i.i.d. $K$-dimensional random vectors; the probability that
$\phi^{0}(j)=e^{\ell}$ is $p_{\ell}$. For $k=0,\ldots,K$, define
the routing
process
\[
\Phi^{k}(N)=\sum_{j=1}^{N}\phi^{k}(j) \qquad\mbox{for }N\in\mathbb{N}.
\]
For each $k=1,\ldots,K$, let $S_{k}$ be a Poisson process with rate
$\nu_{k}$,
and let $G$ be a Poisson process with rate $\alpha$. We assume that
%
%
\begin{equation} \label{eq:independence}\qquad
X^{n}(0),E^{n},S_{1},\ldots,S_{K},\Phi^{0},\ldots,\Phi^{K} \mbox{
and }
G\mbox{ are mutually independent}.
\end{equation}

Let $T_{k}^{n}(t)$ be the cumulative amount of service effort received by
customers in phase $k$ service in $(0,t]$, $B^{n}(t)$ be the cumulative number
of customers who have entered service in $(0,t]$ and $D^{n}(t)$ be the
cumulative number of customers who have completed service in $(0,t]$.
Clearly,
%
%
\begin{equation}\label{eq:Tk}%
T_{k}^{n}(t)=\int_{0}^{t}Z_{k}^{n}(s) \,ds \qquad\mbox{for }t\geq0.
\end{equation}
Then $S_{k}(T_{k}^{n}(t))$ is the cumulative number of phase $k$ service
completions by time~$t$. Also $G(\int_{0}^{t}(X^{n}(s))^{+} \,ds)$ is
the cumulative number of customers who have abandoned the system by
time $t$.
One can check that for $t\geq0$, the processes $X^{n}$ and $Z^{n}$
satisfy the
following dynamical equations:
%
%
\begin{eqnarray}
\label{eq:Zdyn}
Z^{n}(t) & = & Z^{n}(0)+\Phi^{0}(B^{n}(t))+\sum_{k=1}^{K}\Phi^{k}(S_{k}%
(T_{k}^{n}(t)))-S(T^{n}(t)),\\
\label{eq:Xdyn}
X^{n}(t) & = & X^{n}(0)+E^{n}(t)-D^{n}(t)-G \biggl( \int_{0}^{t}(X^{n}%
(s))^{+} \,ds \biggr) ,\\
\label{eq:D1}
D^{n}(t) & = & \sum_{k=1}^{K} \bigl( S_{k}(T_{k}^{n}(t))-e^{\prime}\Phi
^{k}(S_{k}(T_{k}^{n}(t))) \bigr) \nonumber\\[-8pt]\\[-8pt]
& = &-e^{\prime} \Biggl( \sum_{k=1}^{K}\Phi^{k}(S_{k}(T_{k}^{n}(t)))-S(T^{n}%
(t)) \Biggr) ,\nonumber
\end{eqnarray}
where
\[
S(T^{n}(t))= ( S_{1}(T_{1}^{n}(t)),\ldots,S_{K}(T_{K}^{n}(t)) )
^{\prime}.
\]

\subsection{State-process representation}
\label{sec:preAnalysis}

Define the centered processes
\[
\hat{S}(t)=S(t)-\nu t,\qquad \hat{G}(t)=G(t)-\alpha t,\qquad \hat{\Phi}^{\ell
}(N)=\sum_{j=1}^{N} \bigl( \phi^{\ell}(j)-p^{\ell} \bigr)
\]
for $t\geq0$, $\ell=0,\ldots,K$ and $N\in\mathbb{N}$, where
$p^{0}=p$ and
$p^{k}$ is the $k$th column of $P^{\prime}$ for $k=1,\ldots,K$. Setting
%
%
\begin{equation} \label{eq:M}%
M^{n}(t)=\sum_{k=1}^{K}\hat{\Phi}^{k} ( S_{k}(T_{k}^{n}(t)) )
-(I-P^{\prime})\hat{S} ( T^{n}(t) ) ,
\end{equation}
one then has
\[
\sum_{k=1}^{K}\Phi
^{k}(S_{k}(T_{k}^{n}(t)))-S(T^{n}(t))=M^{n}(t)-R\int
_{0}%
^{t}Z^{n}(s) \,ds,
\]
where $R$ is defined in (\ref{eq:R}). By (\ref{eq:Zdyn}) and (\ref{eq:D1}),
%
%
\begin{eqnarray}
\label{eq:B1}
e^{\prime}Z^{n}(t) & = & e^{\prime}Z^{n}(0)+B^{n}(t)-D^{n}(t),\\
\label{eq:DnewRep}%
D^{n}(t) & = & -e^{\prime}M^{n}(t)+e^{\prime}R\int_{0}^{t}Z^{n}(s) \,ds.
\end{eqnarray}
It follows from (\ref{eq:XZ}) and (\ref{eq:Zdyn})--(\ref
{eq:DnewRep}) that
\begin{eqnarray*}
Z^{n}(t) & = & Z^{n}(0)+p(X^{n}(0))^{-}+\hat{\Phi}^{0}(B^{n}(t))-p(X^{n}%
(t))^{-}\\
&&{} +(I-pe^{\prime})M^{n}(t)-(I-pe^{\prime})R\int_{0}^{t}Z^{n}(s) \,ds,\\
X^{n}(t) & = & X^{n}(0)+\hat{E}^{n}(t)+\lambda^{n}t+e^{\prime}M^{n}%
(t)-e^{\prime}R\int_{0}^{t}Z^{n}(s) \,ds\\
&&{} -\hat{G} \biggl( \int_{0}^{t}(X^{n}(s))^{+} \,ds \biggr) -\alpha\int
_{0}^{t}(X^{n}(s))^{+} \,ds.
\end{eqnarray*}
Recall that $\hat{Z}^{n}(t)=Z^{n}(t)-n\gamma$. We then have
\begin{eqnarray*}
\hat{Z}^{n}(t) & = & (I-pe^{\prime})\hat{Z}^{n}(0)+\hat{\Phi}^{0}%
(B^{n}(t))-p(X^{n}(t))^{-}\\
&&{} +(I-pe^{\prime})M^{n}(t)-(I-pe^{\prime})R\int_{0}^{t}\hat{Z}%
^{n}(s) \,ds,\\
X^{n}(t) & = & X^{n}(0)+\hat{E}^{n}(t)+(\lambda^{n}-n\mu)t+e^{\prime}%
M^{n}(t)-e^{\prime}R\int_{0}^{t}\hat{Z}^{n}(s) \,ds\\
&&{} -\hat{G} \biggl( \int_{0}^{t}(X^{n}(s))^{+} \,ds \biggr) -\alpha\int
_{0}^{t}(X^{n}(s))^{+} \,ds,
\end{eqnarray*}
where we have used (\ref{eq:gamma}) and (\ref{eq:R}) in the derivations.
Setting
%
%
\begin{eqnarray}
\label{eq:Un}
U^{n}(t)&=&X^{n}(0)+\hat{E}^{n}(t)+(\lambda^{n}-n\mu)t\nonumber\\[-8pt]\\[-8pt]
&&{}+e^{\prime}%
M^{n}(t)-\hat{G} \biggl( \int_{0}^{t}(X^{n}(s))^{+} \,ds \biggr) ,
\nonumber\\
\label{eq:Vn}%
V^{n}(t)&=&(I-pe^{\prime})\hat{Z}^{n}(0)+\hat{\Phi}^{0}(B^{n}%
(t))+(I-pe^{\prime})M^{n}(t)
\end{eqnarray}
for $t\geq0$, we finally have
%
%
\begin{eqnarray}
\label{eq:rep1}
X^{n}(t) & = & U^{n}(t)-\alpha\int_{0}^{t}(X^{n}(s))^{+} \,ds-e^{\prime}%
R\int_{0}^{t}\hat{Z}^{n}(s) \,ds,\\
\label{eq:rep2}%
\hat{Z}^{n}(t) & = & V^{n}(t)-p(X^{n}(t))^{-}-(I-pe^{\prime})R\int
_{0}^{t}%
\hat{Z}^{n}(s) \,ds.
\end{eqnarray}
By Lemma \ref{lem:Phi}, we have obtained the following representation
for the
state processes:
%
%
\begin{equation} \label{eq:PhiPres}%
(X^{n},\hat{Z}^{n})=\Phi(U^{n},V^{n}).
\end{equation}

\section{Proofs for $G/\PH/n+M$ queues}
\label{sec:exponentialPatience}

This section provides proofs for Theorem
\ref{thm:overloaded} and a special version of Theorem \ref
{thm:critical} when
the patience time distribution is exponential. Section \ref{sec:fluid} first
establishes a fluid limit, which is needed in applying the random-time-change
theorem to prove the theorems in Section \ref{sec:fcl}.

\subsection{Fluid limits}
\label{sec:fluid}

For $t\geq0$, define the fluid-scaled processes $\bar{B}^{n}$, $\bar{D}^{n}$,
$\bar{E}^{n}$, $\bar{T}^{n}$, $\bar{X}^{n}$ and $\bar{Z}^{n}$ via
\begin{eqnarray*}
\bar{B}^{n}(t)&=&\frac{1}{n}B^{n}(t),\qquad \bar{D}^{n}(t)=\frac{1}{n}%
D^{n}(t), \qquad\bar{E}^{n}(t)=\frac{1}{n}E^{n}(t),\\
\bar{T}^{n}(t)&=&\frac{1}{n}T^{n}(t),\qquad \bar{X}^{n}(t)=\frac{1}{n}%
X^{n}(t),\qquad \bar{Z}^{n}(t)=\frac{1}{n}Z^{n}(t).
\end{eqnarray*}

\begin{theorem}
\label{thm:fluid} Consider a sequence of $G/\PH/n+M$ queues satisfying
(\ref{eq:lambdan}) and (\ref{eq:Eclt}). Assume (\ref{eq:OLinitial}) holds.
Then
%
%
\begin{equation} \label{eq:Fluid}%
(\bar{B}^{n},\bar{D}^{n},\bar{E}^{n},\bar{T}^{n},\bar{X}^{n},\bar
{Z}%
^{n})\Rightarrow(\bar{B},\bar{D},\bar{E},\bar{T},\bar{X},\bar
{Z}) \qquad\mbox{as }n\rightarrow\infty,
\end{equation}
where $\bar{B}(t)=\mu t$, $\bar{D}(t)=\mu t$, $\bar{E}(t)=\lambda t$,
$\bar
{T}(t)=\gamma t$, $\bar{X}(t)=q$ and $\bar{Z}(t)=\gamma$ for $t\geq0$.
\end{theorem}
\begin{pf}
For $t\geq0$, let
%
%
\begin{eqnarray} \label{25}%
\bar{M}^{n}(t)&=&\frac{1}{n}M^{n}(t),\qquad \bar{U}^{n}(t)=\frac{1}{n}%
U^{n}(t), \nonumber\\[-8pt]\\[-8pt]
\bar{V}^{n}(t)&=&\frac{1}{n}V^{n}(t), \qquad\bar{L}^{n}(t)=\frac
{1}{n}\hat{Z}^{n}(t).\nonumber
\end{eqnarray}
By (\ref{eq:PhiPres}) and the positively homogeneous property of $\Phi
$, we
have
\[
(\bar{X}^{n},\bar{L}^{n})=\Phi(\bar{U}^{n},\bar{V}^{n}).
\]
Setting
%
%
\begin{equation} \label{eq:bU}%
\bar{U}(t)=q+(\lambda-\mu)t \quad\mbox{and}\quad \bar{V}(t)=0 \qquad\mbox{for
}t\geq0,
\end{equation}
one can check that $\Phi(\bar{U},\bar{V})=(\bar{X},0)$. We are
going to show
that
%
%
\begin{equation} \label{eq:UVfluid}%
(\bar{M}^{n},\bar{U}^{n},\bar{V}^{n})\Rightarrow(0,\bar{U},0)
\qquad\mbox{as
}n\rightarrow\infty.
\end{equation}

Assuming (\ref{eq:UVfluid}), we now complete the proof of the theorem. The
continuity of the map $\Phi$ implies that
\[
(\bar{X}^{n},\bar{L}^{n})=\Phi(\bar{U}^{n},\bar{V}^{n})\quad\Rightarrow\quad
\Phi
(\bar
{U},\bar{V})=(\bar{X},0) \qquad\mbox{as }n\rightarrow\infty.
\]
Since $\bar{Z}^{n}(t)=\bar{L}^{n}(t)+\gamma$ for $t\geq0$, then
$\bar{Z}
^{n}\Rightarrow\bar{Z}$ as $n\rightarrow\infty$, from which one has
$\bar
{T}^{n}\Rightarrow\bar{T}$ as $n\rightarrow\infty$. By (\ref{eq:DnewRep}),
\[
\bar{D}^{n}(t)=-e^{\prime}\bar{M}^{n}(t)+e^{\prime}R\int
_{0}^{t}\bar{Z}%
^{n}(s) \,ds.
\]
Since $e^{\prime}R\int_{0}^{t}\bar{Z}(s) \,ds=\mu t$ for $t\geq0$, by the
continuous-mapping theorem $\bar{D}^{n}\Rightarrow\bar{D}$ as
$n\rightarrow
\infty$. The convergence of $\bar{D}^{n}$ and (\ref{eq:B1}) imply that
$\bar{B}^{n}\Rightarrow\bar{B}$ as $n\rightarrow\infty$, and $\bar{B}$
satisfies
\[
e^{\prime}\bar{Z}(t)=e^{\prime}\bar{Z}(0)+\bar{B}(t)-\bar{D}(t)
\qquad\mbox{for }t\geq0.
\]
Since $\bar{Z}(t)=\bar{Z}(0)=\gamma$, we conclude that $\bar
{B}(t)=\mu
t$ for
$t\geq0$. By assumptions (\ref{eq:lambdan}) and (\ref{eq:Eclt}), for each
$T>0$,
%
%
\begin{equation}\label{eq:SLLNE}%
{\frac{1}{n}\sup_{0\leq t\leq T} }\vert\hat{E}^{n}(t) \vert
\Rightarrow0 \qquad\mbox{as }n\rightarrow\infty,
\end{equation}
which implies that $\bar{E}^{n}\Rightarrow\bar{E}$. This proves the theorem
when (\ref{eq:UVfluid}) holds.

It remains to prove (\ref{eq:UVfluid}). By the functional strong law
of large
numbers \mbox{(FSLLN)},
%
%
\begin{eqnarray} \label{eq:SLLNGSPhi}%
{\frac{1}{n}\sup_{0\leq t\leq T}} \vert\hat{G}(nt) \vert
&\Rightarrow& 0,\qquad {\frac{1}{n}\sup_{0\leq t\leq T} }\vert\hat
{S}(nt) \vert\Rightarrow0, \nonumber\\[-8pt]\\[-8pt]
{\frac{1}{n}\sup_{0\leq t\leq T} }\vert
\hat{\Phi}^{k}(\lfloor nt\rfloor) \vert&\Rightarrow&0\nonumber
\end{eqnarray}
as $n\rightarrow\infty$, for $k=0,\ldots,K$. Let $\bar
{S}^{n}(t)=S^{n}(nt)/n$
for $t\geq0$. The FSLLN also leads to $\bar{S}^{n}\Rightarrow\bar
{S}$ as
$n\rightarrow\infty$ where $\bar{S}(t)=\nu t$ for $t\geq0$. This, together
with (\ref{eq:M}), (\ref{eq:SLLNGSPhi}) and the fact $\bar
{T}_{k}^{n}(t)\leq
t$ for $t\geq0$, implies that
%
%
\begin{equation} \label{eq:Mtozero}%
{\sup_{0\leq t\leq T}} \vert\bar{M}^{n}(t) \vert\Rightarrow
0 \qquad\mbox{as }n\rightarrow\infty.
\end{equation}
Note that $\bar{B}^{n}(t)\leq(\bar{X}^{n}(0))^{+}+\bar{E}^{n}(t)$. By
(\ref{eq:Zinitial}) and the convergence of $\bar{E}^{n}$, the
sequence of
processes $\{\bar{B}^{n},n\in\mathbb{N}\}$ is stochastically bounded,
that is,
for each \mbox{$T>0$},
\[
\lim_{a\rightarrow\infty}\limsup_{n\rightarrow\infty}\mathbb{P} \Bigl[
\sup_{0\leq t\leq T}\bar{B}^{n}(t)>a \Bigr] =0.
\]
Using this and (\ref{eq:SLLNGSPhi}), we deduce that%
%
%
\begin{equation} \label{eq:Phi0Zero}
\sup_{0\leq t\leq T}\frac{1}{n}\hat{\Phi}^{0} ( B^{n}(t) )
\Rightarrow0 \qquad\mbox{as }n\rightarrow\infty.%
\end{equation}
Condition (\ref{eq:Zinitial}) implies that $\hat
{Z}^{n}(0)/n\Rightarrow
0$ as
$n\rightarrow\infty$, which, together with (\ref{eq:Mtozero}) and
(\ref{eq:Phi0Zero}), leads to $\bar{V}^{n}\Rightarrow0$ as
$n\rightarrow
\infty$. Since $\sup_{0\leq t\leq T}(\bar{X}^{n}(t))^{+}\leq(\bar
{X}%
^{n}(0))^{+}+\bar{E}^{n}(T)$, one can argue similarly that $\bar{U}%
^{n}\Rightarrow\bar{U}$ as $n\rightarrow\infty$. Hence, we have shown
(\ref{eq:UVfluid}) holds and thus proved the theorem.
\end{pf}

\subsection{Diffusion limits}
\label{sec:fcl}

In this section, we prove Theorems \ref{thm:critical} and \ref
{thm:overloaded}%
, under the assumption that the patience times are exponentially
distributed.\vspace*{1pt}

Define the diffusion-scaled processes $\tilde{G}^{n}$, $\tilde
{S}^{n}$ and
$\tilde{\Phi}^{0},\ldots,\tilde{\Phi}^{K}$ via%
\[
\tilde{G}^{n}(t)=\frac{1}{\sqrt{n}}\hat{G}(nt),\qquad \tilde
{S}^{n}(t)=\frac
{1}{\sqrt{n}}\hat{S}(nt),\qquad \tilde{\Phi}^{k,n}(t)=\frac{1}{\sqrt{n}}
\hat{\Phi}^{k}(\lfloor nt\rfloor)
\]
for $t\geq0$ and $k=0,\ldots,K$. By the FCLT, one has
\[
(\tilde{G}^{n},\tilde{S}^{n},\tilde{\Phi}^{0,n},\ldots,\tilde
{\Phi}%
^{K,n})\Rightarrow(\tilde{G},\tilde{S},\tilde{\Phi}^{0},\ldots
,\tilde
{\Phi
}^{K}) \qquad\mbox{as }n\rightarrow\infty,
\]
where $\tilde{G}$ is a one-dimensional driftless Brownian motion, and
$\tilde{S}$ and $\tilde{\Phi}^{0},\ldots,\tilde{\Phi}^{K}$ are
$K$-dimensional
driftless Brownian motions. As mentioned previously, the variance of
$\tilde{G}$ is $\alpha$, the covariance matrix for $\tilde{S}$ is
$\operatorname{diag}
(\nu
)$, and for $k=0,\ldots,K$, the covariance matrix for $\tilde{\Phi
}^{k}$ is
$H^{k}$ given by (\ref{eq:Hk}). By the FCLT assumption (\ref{eq:Eclt})
for the
arrival process $E^{n}$, the initial condition (\ref{eq:OLinitial}),
and the
independence assumption (\ref{eq:independence}), one has
%
%
\begin{eqnarray} \label{eq:FCLT}%
&&(\tilde{X}^{n}(0),\tilde{Z}^{n}(0),\tilde{E}^{n},\tilde
{G}^{n},\tilde{S}
^{n},\tilde{\Phi}^{0,n},\ldots,\tilde{\Phi}^{K,n})\nonumber\\[-8pt]\\[-8pt]
&&\qquad\Rightarrow
(\tilde
{X}(0),\tilde{Z}(0),\tilde{E},\tilde{G},\tilde{S},\tilde{\Phi}^{0}
,\ldots,\tilde{\Phi}^{K})\nonumber
\end{eqnarray}
as $n\rightarrow\infty$. The components of $(\tilde{E},\tilde
{G},\tilde
{S},\tilde{\Phi}^{0},\ldots,\tilde{\Phi}^{K})$ are mutually
independent, and
they are independent of $(\tilde{X}(0),\tilde{Z}(0))$.\vspace*{1pt}

Let $\hat{U}^{n}(t)=U^{n}(t)-n\bar{U}(t)$, and define the diffusion-scaled
processes
\[
\tilde{U}^{n}(t)=\frac{1}{\sqrt{n}}\hat{U}^{n}(t) \quad\mbox{and}\quad
\tilde{V}^{n}(t)=\frac{1}{\sqrt{n}}V^{n}(t) \qquad\mbox{for }t\geq0,
\]
where $\bar{U}$ is defined in (\ref{eq:bU}). We now have the
following lemma.
\begin{lemma}
\label{lem:BMdiffusion} Consider a sequence of $G/\PH/n+M$ queues satisfying
(\ref{eq:lambdan}) and (\ref{eq:Eclt}). Assume that (\ref{eq:OLinitial})
holds. Then
\[
(\tilde{U}^{n},\tilde{V}^{n})\Rightarrow(\tilde{U},\tilde{V})
\qquad\mbox{as }n\rightarrow\infty,
\]
where $(\tilde{U},\tilde{V})$ is a $(K+1)$-dimensional Brownian
motion defined
by (\ref{eq:tUOvr}) and (\ref{eq:tV}).
\end{lemma}
\begin{pf}
By (\ref{eq:Un}) and (\ref{eq:Vn}),
%
%
\begin{eqnarray}
\label{eq:tilUn}
\tilde{U}^{n}(t)&=&\tilde{X}^{n}(0)+\tilde{E}^{n}(t)+\sqrt{n} \biggl(
\frac
{1}{n}\lambda^{n}-\lambda\biggr)\nonumber\\[-8pt]\\[-8pt]
&&{} +e^{\prime}\tilde{M}^{n}(t)-\tilde{G}%
^{n} \biggl( \int_{0}^{t}(\bar{X}^{n}(s))^{+} \,ds \biggr),\nonumber\\
\label{eq:tilVn}%
\tilde{V}^{n}(t)&=&(I-pe^{\prime})\tilde{Z}^{n}(0)+\tilde{\Phi
}^{0,n}(\bar
{B}^{n}(t))+(I-pe^{\prime})\tilde{M}^{n}(t),
\end{eqnarray}
where
\[
\tilde{M}^{n}(t)=\frac{1}{\sqrt{n}}M^{n}(t)=\sum_{k=1}^{K}\tilde
{\Phi}%
^{k,n}(\bar{S}_{k}^{n}(\bar{T}_{k}^{n}(t)))-(I-P^{\prime})\tilde
{S}^{n}%
(\bar{T}^{n}(t))
\]
and $\bar{S}^{n}(t)=S(nt)/n$ for $t\geq0$. By the FSLLN, $\bar{S}%
^{n}\Rightarrow\bar{S}$ as $n\rightarrow\infty$, where $\bar
{S}(t)=\nu
t$ for
$t\geq0$. The lemma now follows from (\ref{eq:FCLT}), Theorem \ref
{thm:fluid},
the continuous-mapping theorem and the random-time-change theorem.
\end{pf}
\begin{pf*}{Proof of Theorem \ref{thm:critical} \textup{(Assuming an
exponential patience
time distribution)}}
Since\vadjust{\goodbreak} $\rho=1$, it follows that $q=0$ and
$\lambda=\mu$. Then $\bar{U}(t)=0$ for $t\geq0$. It follows from
the state-process representation (\ref{eq:PhiPres}) and the
positively homogeneous property of the map $\Phi$ that
\[
(\tilde{X}^{n},\tilde{Z}^{n})=\Phi(\tilde{U}^{n},\tilde{V}^{n}).
\]
The theorem now follows from Lemma \ref{lem:BMdiffusion} and the
continuous-mapping theorem.
\end{pf*}

Although condition (\ref{eq:initial}) is not explicitly required
in the above proof, it can be deduced by using initial condition
(\ref{eq:Zinitial}) and the assumption that the patience times of
all customers, including those in queue initially, are i.i.d. following
an exponential distribution.

Before proving Theorem \ref{thm:overloaded}, we first establish a
lemma. It
says that when $\rho>1$, the number of idle servers goes to zero under
diffusion scaling.
\begin{lemma}
\label{lem:zeroIdle} Let $I^{n}(t)=(X^{n}(t))^{-}$ and $\tilde{I}%
^{n}(t)=I^{n}(t)/\sqrt{n}$ for $t\geq0$. Then under the conditions of
Theorem \ref{thm:overloaded},
\[
\tilde{I}^{n}\Rightarrow0 \qquad\mbox{as }n\rightarrow\infty.
\]
\end{lemma}
\begin{pf}
It follows from (\ref{eq:rep1}) and (\ref{eq:rep2}) that
\begin{eqnarray*}
\frac{1}{\sqrt{n}}X^{n}(t) & = & \sqrt{n}\bar{U}(t)+\tilde{U}^{n}%
(t)-\frac{\alpha}{\sqrt{n}}\int_{0}^{t}(X^{n}(s))^{+} \,ds-e^{\prime
}R\int
_{0}^{t}\tilde{Z}^{n}(s) \,ds,\\
\tilde{Z}^{n}(t) & = & \tilde{V}^{n}(t)-\frac{p}{\sqrt
{n}}(X^{n}(t))^{-}%
-(I-pe^{\prime})R\int_{0}^{t}\tilde{Z}^{n}(s) \,ds.
\end{eqnarray*}
Therefore, by Lemma \ref{lem:Phi}
\[
\biggl( \frac{1}{\sqrt{n}}X^{n},\tilde{Z}^{n} \biggr) =\Phi\bigl(\tilde{U}%
^{n}+\sqrt{n}\bar{U},\tilde{V}^{n}\bigr).
\]
By the Lipschitz continuity property (\ref{eq:LipschitzPhi}) of the map
$\Phi
$, for any $T>0$, there exists a constant $C_{T}^{1}>0$ such that
\[
\sup_{0\leq t\leq T} \bigl\vert\Phi\bigl(\tilde{U}^{n}+\sqrt{n}\bar
{U},\tilde
{V}^{n}\bigr)(t)-\Phi\bigl(\sqrt{n}\bar{U},0\bigr)(t) \bigr\vert\leq C_{T}^{1}\sup
_{0\leq
t\leq
T} \{ |\tilde{U}^{n}(t)|+|\tilde{V}^{n}(t)| \}
\]
for all $n$ and all sample paths. One can check that $\Phi(\sqrt
{n}\bar
{U},0)=(\sqrt{n}q,0)$. Therefore,
%
%
\begin{equation}\label{eq:Xlowerbound}%
\inf_{0\leq t\leq T}\frac{1}{\sqrt{n}}X^{n}(t)\geq\sqrt
{n}q-C_{T}^{1}\sup_{0\leq
t\leq T} \{ |\tilde{U}^{n}(t)|+|\tilde{V}^{n}(t)| \} .
\end{equation}
By Lemma \ref{lem:BMdiffusion},
\[
\lim_{a\rightarrow\infty}\limsup_{n\rightarrow\infty}\mathbb{P} \Bigl[
\sup_{0\leq t\leq T} \{ |\tilde{U}^{n}(t)|+|\tilde{V}^{n}(t)| \}
>a \Bigr] =0,
\]
which, together with (\ref{eq:Xlowerbound}), implies that $\sup
_{0\leq
t\leq
T}\tilde{I}^{n}(t)\Rightarrow0$ as $n\rightarrow\infty$. The lemma
is proved.
\end{pf}
\begin{pf*}{Proof of Theorem \ref{thm:overloaded}}
It follows from (\ref{eq:rep1}) and
(\ref{eq:rep2}) that
\begin{eqnarray*}
\tilde{X}^{n}(t) & = & \tilde{U}^{n}(t)-\frac{\alpha}{\sqrt{n}}\int
_{0}%
^{t}(X^{n}(s))^{-} \,ds-\alpha\int_{0}^{t}\tilde{X}^{n}(s)
\,ds-e^{\prime}%
R\int_{0}^{t}\tilde{Z}^{n}(s) \,ds,\\
\tilde{Z}^{n}(t) & = & \tilde{V}^{n}(t)-\frac{p}{\sqrt
{n}}(X^{n}(t))^{-}%
-(I-pe^{\prime})R\int_{0}^{t}\tilde{Z}^{n}(s) \,ds.
\end{eqnarray*}
Let
\[
\delta^{n}(t)=\frac{\alpha}{\sqrt{n}}\int_{0}^{t}(X^{n}(s))^{-} \,ds
\quad\mbox{and}\quad \epsilon^{n}(t)=\frac{p}{\sqrt{n}}(X^{n}(t))^{-}\qquad
\mbox{for }t\geq0.
\]
By Lemma \ref{lem:Psi},
\[
(\tilde{X}^{n},\tilde{Z}^{n})=\Psi(\tilde{U}^{n}-\delta^{n},\tilde
{V}%
^{n}-\epsilon^{n}).
\]
By Lemma \ref{lem:zeroIdle},
%
%
\begin{equation}\label{eq:Idletozero}%
(\delta^{n},\epsilon^{n})\Rightarrow(0,0) \qquad\mbox{as }n\rightarrow
\infty.
\end{equation}
The theorem follows from Lemma \ref{lem:BMdiffusion}, (\ref{eq:Idletozero})
and the continuity of the map $\Psi$.
\end{pf*}

\section{Proofs for critically loaded $G/\PH/n+\GI$ queues}
\label{sec:criGI}

In this section, we prove Theorem \ref{thm:critical}
for a
general patience time distribution. Consider a sequence of $G/\PH/n+\GI$ queues
indexed by $n$. Our starting point is the perturbed system described in
Section \ref{sec:perturbed} with the following modification: each
customer in
queue can abandon the system, not just the leading customer; when
a customer's waiting time in the real FIFO queue exceeds his patience time,
the customer abandons the system. By the same argument as in Section
\ref{sec:perturbed}, for each $n$, the modified perturbed system is equivalent
in distribution to the original $G/\PH/n+\GI$ queue. In particular, the system
equations (\ref{eq:Zdyn})--(\ref{eq:D1}) derived in Section \ref
{sec:dynamics}
hold, except that (\ref{eq:Xdyn}) is modified as follows:
%
%
\begin{equation} \label{eq:GIXdyn}%
X^{n}(t)=X^{n}(0)+E^{n}(t)-D^{n}(t)-A^{n}(t),
\end{equation}
where $A^{n}(t)$ denotes the cumulative number of customers that have
abandoned the system by time $t$. We call $A^{n}=\{A^{n}(t),t\geq0\}$ the
\textit{abandonment-count process} in the $n$th system.

With systems equations (\ref{eq:Zdyn}), (\ref{eq:GIXdyn}) and (\ref{eq:D1}),
one can derive representation (\ref{eq:PhiPres})
\[
(X^{n}, \hat Z^{n}) = \Phi(U^{n}, V^{n})
\]
with $U^{n}$ modified as
%
%
\begin{eqnarray} \label{eq:GIUn}%
U^{n}(t)&=&X^{n}(0)+\hat{E}^{n}(t)+(\lambda^{n}-n\mu)t+e^{\prime}M^{n}
(t)\nonumber\\[-8pt]\\[-8pt]
&&{}-A^{n}(t)+\alpha\int_{0}^{t}(X^{n}(s))^{+} \,ds.\nonumber
\end{eqnarray}
The derivation is identical to the one in Section \ref{sec:preAnalysis}
and is
not repeated here.\vadjust{\goodbreak}

Before we prove Theorem \ref{thm:critical}, we state two lemmas, which
will be
proved at the end of this section. The first lemma follows from a main result
in \citet{DaiHe09}. In the lemma, the diffusion-scaled abandonment-count
process $\tilde{A}^{n}=\{\tilde{A}^{n}(t),t\geq0\}$ is defined by
%
%
\begin{equation}\label{eq:tilAn}%
\tilde{A}^{n}(t)=\frac{1}{\sqrt{n}}A^{n}(t) \qquad\mbox{for }t\geq0.
\end{equation}

\begin{lemma}
\label{lem:DaiHe} Under the conditions of Theorem \ref{thm:critical},
for any $T>0$,
%
%
\begin{equation}\label{eq:DaiHe09}%
{\sup_{0\leq t\leq T}}\biggl \vert\tilde{A}^{n}(t)-\alpha\int
_{0}^{t}(\tilde
{X}^{n}(s))^{+} \,ds \biggr\vert\Rightarrow0 \qquad\mbox{as }n\rightarrow\infty.
\end{equation}
\end{lemma}

The next lemma is a generalization of the fluid limit theorem (Theorem
\ref{thm:fluid}) to general patience time distributions, but with the
restriction that $\rho=1$.

\begin{lemma}
\label{lem:fluid+GI} Under the conditions of Theorem \ref{thm:critical},
%
%
\begin{equation} \label{eq:Fluid+GI}%
(\bar{B}^{n},\bar{D}^{n},\bar{E}^{n},\bar{T}^{n},\bar{X}^{n},\bar
{Z}%
^{n})\Rightarrow(\bar{B},\bar{D},\bar{E},\bar{T},\bar{X},\bar
{Z}) \qquad\mbox{as }n\rightarrow\infty,
\end{equation}
where $\bar{B}^{n}$, $\bar{D}^{n}$, $\bar{E}^{n}$, $\bar{T}^{n}$,
$\bar
{X}%
^{n}$ and $\bar{Z}^{n}$ are fluid-scaled processes defined at the beginning
of Section \ref{sec:fluid}, and $\bar{B}(t)=\mu t$, $\bar{D}(t)=\mu t$,
$\bar{E}(t)=\lambda t$, $\bar{T}(t)=\gamma t$, $\bar{X}(t)=0$ and
$\bar
{Z}(t)=\gamma$ for $t\geq0$.
\end{lemma}

\begin{pf*}{Proof of Theorem \ref{thm:critical}}
Using the representation
(\ref{eq:PhiPres}) with $U^{n}$ given by (\ref{eq:GIUn}), one has
\[
(\tilde{X}^{n},\tilde{Z}^{n})=\Phi(\tilde{U}^{n},\tilde{V}^{n}),
\]
where%
%
%
\begin{eqnarray}\label{eq:tUnCri}%
\tilde{U}^{n}(t)&=&\tilde{X}^{n}(0)+\tilde{E}^{n}(t)+\sqrt{n} \biggl( \frac
{1}%
{n}\lambda^{n}-\lambda\biggr)
+e^{\prime}\tilde{M}^{n}(t)\nonumber\\[-8pt]\\[-8pt]
&&{} - \biggl( \tilde
{A}^{n}(t)-\alpha\int_{0}^{t}(\tilde{X}^{n}(s))^{+} \,ds \biggr)\nonumber
\end{eqnarray}
and $\tilde{V}^{n}$ is given by (\ref{eq:tilVn}). By Lemma \ref
{lem:Phi}, the
map $\Phi$ is continuous. Thus, to prove the theorem, it suffices to prove
that
%
%
\begin{equation}\label{eq:tUnVnConverge}%
(\tilde{U}^{n},\tilde{V}^{n})\Rightarrow(\tilde{U},\tilde{V}),
\end{equation}
where $(\tilde{U},\tilde{V})$ is the $(K+1)$-dimensional Brownian motion
defined by (\ref{eq:tU}) and (\ref{eq:tV}). The convergence
(\ref{eq:tUnVnConverge}) follows from the proof of Lemma \ref{lem:BMdiffusion}
with the following two modifications. First, the last term of $\tilde{U}^{n}$
in (\ref{eq:tUnCri}) is
\[
\biggl( \tilde{A}^{n}(t)-\alpha\int_{0}^{t}(\tilde{X}^{n}(s))^{+} \,ds \biggr)
\]
instead of
\[
\tilde{G}^{n} \biggl( \int_{0}^{t}(\bar{X}^{n}(s))^{+} \,ds
\biggr)\vadjust{\goodbreak}
\]
in (\ref{eq:tilUn}). We apply Lemma \ref{lem:DaiHe} to conclude that
the last
term in (\ref{eq:tUnCri}) converges to zero in distribution. Second,
we use
(\ref{eq:Fluid+GI}) instead of (\ref{eq:Fluid}) in order to apply the
random-time-change theorem to finish the proof of (\ref{eq:tUnVnConverge}).
\end{pf*}
\begin{remark*}
It follows immediately from Theorem \ref{thm:critical} and Lemma
\ref{lem:DaiHe} that under the conditions of Theorem \ref{thm:critical},
%
%
\begin{equation}\label{eq:AbConverge}%
\tilde{A}^{n}\Rightarrow\tilde{A} \qquad\mbox{as }n\rightarrow\infty,
\end{equation}
where
\[
\tilde{A}(t)=\alpha\int_{0}^{t}(\tilde{X}(s))^{+} \,ds \qquad\mbox{for
}t\ge0.
\]
\end{remark*}
\begin{pf*}{Proof of Lemma \ref{lem:DaiHe}}
We use Theorem 2.1 of \citet{DaiHe09} to prove
the lemma. In order to apply the theorem, we need only verify that the
sequence of
diffusion-scaled queue-length processes is stochastically bounded, that is,
for any $T>0$,%
%
%
\begin{equation}\label{eq:GIsb}%
\lim_{a\rightarrow\infty}\limsup_{n\rightarrow\infty}\mathbb{P} \biggl[
\sup_{0\leq t\leq T}\frac{1}{\sqrt{n}}(X^{n}(t))^{+}>a \biggr] =0.
\end{equation}
Theorem 2.2 of \citet{DaiHe09} states a comparison result: the queue
length at
any time in a $G/G/n+G$ queue is dominated by the queue length in the
corresponding $G/G/n$ queue without abandonment. Thus, (\ref{eq:GIsb}) is
implied by the stochastic boundedness of the sequence of diffusion-scaled
queue-length processes in the corresponding $G/\PH/n$ queues. Examining the
proof of Theorem \ref{thm:critical} in \mbox{Section \ref{sec:fcl}} for an
exponential patience time distribution with rate $\alpha>0$, one concludes
that Theorem \ref{thm:critical} holds for the corresponding $G/\PH/n$ queues
without abandonment by setting $\alpha=0$. As a consequence, the
sequence of
diffusion-scaled queue-length processes in the $G/\PH/n$ queues is
stochastically bounded.
\end{pf*}
\begin{pf*}{Proof of Lemma \ref{lem:fluid+GI}}
The proof of the lemma follows the proof of
Theorem \ref{thm:fluid} with the following two modifications. First,
$\bar{U}$
in (\ref{eq:bU}) becomes zero in the current case because $\rho=1$. Second,
$U^{n}$ has the representation (\ref{eq:tUnCri}) instead of (\ref
{eq:Un}). In
Theorem \ref{thm:fluid}, to prove $\bar{U}^{n}\Rightarrow0$ as
$n\rightarrow
\infty$ we used the fact
\[
{\frac{1}{n}\sup_{0\leq t\leq T}}|\hat{G}(nt)|\Rightarrow0 \qquad\mbox{as
}n\rightarrow\infty,
\]
which is proved in (\ref{eq:SLLNGSPhi}). Here, we need%
\[
\frac{1}{n}\sup_{0\leq t\leq T} \biggl\vert A^{n}(t)-\alpha\int_{0}^{t}%
(X^{n}(s))^{+} \,ds \biggr\vert\Rightarrow0 \qquad\mbox{as }n\rightarrow
\infty,
\]
which holds because of Lemma \ref{lem:DaiHe}.\vspace*{-14pt}
\end{pf*}

\begin{appendix}
\section{A continuous map}
\label{sec:Map}

Let $K\in\mathbb{N}$ be a fixed positive integer. Given functions
$h_{1}\dvtx\mathbb{R}%
^{K+1}\rightarrow\mathbb{R}$, $h_{2}\dvtx\mathbb{R}^{K+1}\rightarrow
\mathbb
{R}%
^{K}$ and\vadjust{\goodbreak} $g\dvtx\mathbb{R}\rightarrow\mathbb{R}^K$, we wish to define a map
$\Upsilon\dvtx\mathbb{D}^{K+1}\rightarrow\mathbb{D}^{K+1}$. For each
$y=(y_{1},y_{2})\in\mathbb{D}^{K+1}$
with $y_{1}(t)\in\mathbb{R}$ and $y_{2}(t)\in\mathbb{R}^{K}$ for
$t\geq
0$, $\Upsilon(y)$
is defined to be any $x=(x_{1},x_{2})\in\mathbb{D}^{K+1}$ with
$x_{1}(t)\in\mathbb{R}$
and $x_{2}(t)\in\mathbb{R}^{K}$ for $t\geq0$ that satisfies
%
%
\begin{eqnarray}
\label{eq:mapPhi21}
x_{1}(t)&=&y_{1}(t)+\int_{0}^{t}h_{1}(x(s)) \,ds,\\
\label{eq:mapPhi22}%
x_{2}(t)&=&y_{2}(t)+\int_{0}^{t}h_{2}(x(s)) \,ds+g(x_{1}(t))
\end{eqnarray}
for $t\geq0$. We assume that $h_{1}$, $h_{2}$ and $g$ are Lipschitz
continuous. For a function $f\dvtx\mathbb{R}^{d}\rightarrow\mathbb
{R}^{m}$ with
$d,m\in\mathbb{N}$, it is said to be \textit{Lipschitz continuous} with
Lipschitz constant $c>0$ if
\[
|f(u)-f(v)|\leq c|u-v| \qquad\mbox{for }u,v\in\mathbb{R}^{d}
\]
(recall that $|u|={\max_{1\leq k \leq d}}|u_{k}|$ denotes the maximum norm
of $u$). The function $f$ is said to be \textit{positively homogeneous} if
\[
f(au)=af(u) \qquad\mbox{for any }a>0\mbox{ and }u\in\mathbb{R}^{d}.
\]
Given $d\in\mathbb{N}$, $x\in\mathbb{D}^{d}$ and $T>0$, set $\Vert
x\Vert_{T}%
={\sup_{0\leq t\leq T}}\vert x(t)\vert$.

The following lemma establishes the existence and the continuity of the map
$\Upsilon$.
\begin{lemma}
\label{lem:maps} Assume that $h_{1}$, $h_{2}$ and $g$ are Lipschitz
continuous. \textup{(a)} For each $y=(y_{1},y_{2})\in\mathbb{D}^{K+1}$
with $y_{1}(t)\in\mathbb{R}$ and $y_{2}(t)\in\mathbb{R}^{K}$ for
$t\geq0$,
there exists a unique $x=(x_{1},x_{2})\in\mathbb{D}^{K+1}$ with
$x_{1}(t)\in\mathbb{R}$ and $x_{2}(t)\in\mathbb{R}^{K}$ for $t\geq
0$ that
satisfies (\ref{eq:mapPhi21}) and (\ref{eq:mapPhi22}). \textup{(b)} The map
$\Upsilon\dvtx
\mathbb{D}^{K+1}\rightarrow\mathbb{D}^{K+1}$ is
Lipschitz continuous in the sense that for each $T>0$, there exists a constant
$C_{T}>0$ such that
\[
\Vert\Upsilon(y)-\Upsilon(\tilde{y})\Vert_{T}\leq C_{T}\Vert
y-\tilde{y}
\Vert_{T} \qquad\mbox{for any }y,\tilde{y}\in\mathbb{D}^{K+1}.
\]
\textup{(c)} The map $\Upsilon$ is continuous when the domain $\mathbb{D}^{K+1}$
and the range $\mathbb{D}^{K+1}$ are both
endowed with the Skorohod $J_{1}$-topology. \textup{(d)} If, in addition,
$h_{1}$, $h_{2}$ and $g$ are assumed to be positively homogeneous, then the
map $\Upsilon$ is positively homogeneous in the sense that
%
\[
\Upsilon(ay)=a\Upsilon(y) \qquad\mbox{for each }a>0\mbox{ and each }%
y\in\mathbb{D}^{K+1}.
\]
\end{lemma}
\begin{pf}
Assume that $h_{1}$, $h_{2}$ and $g$ are Lipschitz continuous with Lipschitz
constant $c>0$. Let $y=(y_{1},y_{2})\in\mathbb{D}^{K+1}$ be given.
Let $T>0$
be fixed for the moment. Define $x^{0}=y$ and for each $n\in\mathbb{Z}_{+}$,
let $x^{n+1}=(x_{1}^{n+1},x_{2}^{n+1})$ be defined via
\begin{eqnarray*}
x_{1}^{n+1}(t)&=&y_{1}(t)+\int_{0}^{t}h_{1}(x^{n}(s)) \,ds,\\
x_{2}^{n+1}(t)&=&y_{2}(t)+\int_{0}^{t}h_{2}(x^{n}(s)) \,ds+g(x_{1}^{n+1}(t))
\end{eqnarray*}
for $t\in\lbrack0,T]$. Setting
\[
X^{(n)}(t)=\Vert x^{n+1}-x^{n}\Vert_{t},\vadjust{\goodbreak}
\]
because
\begin{eqnarray*}
x_{2}^{n+1}(t)-x_{2}^{n}(t) & = & \int_{0}^{t} \bigl( h_{2}(x^{n}(s))-h_{2}%
(x^{n-1}(s)) \bigr) \,ds\\
&&{} +g \biggl( y_{1}(t)+\int_{0}^{t}h_{1}(x^{n}(s)) \,ds \biggr) \\
&&{} -g \biggl(
y_{1}(t)+\int_{0}^{t}h_{1}(x^{n-1}(s)) \,ds \biggr)
\end{eqnarray*}
for $t\in\lbrack0,T]$, one has
\[
X^{(n+1)}(t)\leq(c+c^{2})\int_{0}^{t}X^{(n)}(s) \,ds \qquad\mbox{for }%
t\in\lbrack0,T].
\]
Then, by Lemma 11.3 in \citet{MandelStrong},
\[
X^{n+1}(t)\leq(c+c^{2})\frac{T^{n}}{n!}\sup_{0\leq s\leq t}X^{(0)}%
(s) \qquad\mbox{for }t\in\lbrack0,T].
\]
Therefore, similarly to (11.22) in \citet{MandelStrong}, $\{
x^{n},n\in
\mathbb{N}\}$ is a Cauchy sequence under the uniform norm \mbox{$\Vert\cdot
\Vert
_{T}$}. Since $(\mathbb{D}([0,T],\mathbb{R}^{K+1})$, \mbox{$\Vert\cdot\Vert_{T})$}
is a
complete metric space (being a closed subset of the Banach space of bounded
functions defined from $[0,T]$ into $\mathbb{R}^{K+1}$ and endowed
with the
uniform norm), $\{x^{n},n\in\mathbb{N}\}$ has a limit $x$ that is in
$\mathbb{D}([0,T],\mathbb{R}^{K+1})$. One can check that $x$ satisfies
(\ref{eq:mapPhi21}) and (\ref{eq:mapPhi22}) for $t\in\lbrack0,T]$.
This proves
the existence of the map $\Upsilon$ from $\mathbb{D}([0,T],\mathbb{R}^{K+1})$
to $\mathbb{D}([0,T],\mathbb{R}^{K+1})$.

Now we prove that the map from $\mathbb{D}([0,T],\mathbb{R}^{K+1})$ to
$\mathbb{D}([0,T],\mathbb{R}^{K+1})$ is Lipschitz continuous with
respect to
the uniform norm. Assume that $y,\tilde{y}\in\mathbb
{D}([0,T],\mathbb{R}
^{K+1})$. Let $\Upsilon(y)$ be any solution $x$ such that $x$ and $y$ satisfy
(\ref{eq:mapPhi21}) and (\ref{eq:mapPhi22}) on $[0,T]$. Similarly, let
$\Upsilon(\tilde{y})$ be any solution associated with $\tilde{y}$. Setting
$x=(x_{1},x_{2})=\Upsilon(y)$ and $\tilde{x}=(\tilde{x}_{1},\tilde
{x}%
_{2})=\Upsilon(\tilde{y})$, then for any $t\in\lbrack0,T]$,
\begin{eqnarray*}
|x_{1}(t)-\tilde{x}_{1}(t)| & \leq & |y(t)-\tilde{y}(t)|+c\int_{0}^{t}%
|\Upsilon(y)(s)-\Upsilon(\tilde{y})(s)| \,ds,\\
|x_{2}(t)-\tilde{x}_{2}(t)| & \leq & (1+c)|y(t)-\tilde
{y}(t)|\\
&&{}+(c+c^{2})\int
_{0}^{t}|\Upsilon(y)(s)-\Upsilon(\tilde{y})(s)| \,ds.
\end{eqnarray*}
Hence,
\begin{eqnarray}
&&|\Upsilon(y)(t)-\Upsilon(\tilde{y})(t)|
\nonumber\\
&&\qquad\leq(1+c)|y(t)-\tilde
{y}(t)|
+(c+c^{2}%
)\int_{0}^{t}|\Upsilon(y)(s)-\Upsilon(\tilde{y})(s)| \,ds \nonumber\\
\eqntext{\mbox{for
}%
t\in[0,T].}
\end{eqnarray}
By Corollary 11.2 in \citet{MandelStrong}
\[
\Vert\Upsilon(y)-\Upsilon(\tilde{y})\Vert_{T}\leq(1+c)\Vert
y-\tilde{y}%
\Vert_{T}\exp\bigl((c+c^{2})T\bigr).
\]
Hence, $\Upsilon$ is Lipschitz continuous, which implies part (b) of the
lemma. The Lipschitz continuity of $\Upsilon$ as a map from $\mathbb
{D}%
([0,T],\mathbb{R}^{K+1})$ to $\mathbb{D}([0,T],\mathbb{R}^{K+1})$
shows that
it is well defined on $[0,T]$. Since $T>0$ is arbitrary, $\Upsilon$ as
a map
from $\mathbb{D}^{K+1}$ to $\mathbb{D}^{K+1}$ is well defined. This proves
part (a) of the lemma.

Next we prove the continuity of $\Upsilon$ provided that $\mathbb
{D}^{K+1}$ is
endowed with the Skorohod $J_{1}$-topology [see, e.g.,
Section 3 of \citet{Whitt02}]. Consider a sequence $\{y^{n},n\in
\mathbb{N}\}$
and $y$ in $\mathbb{D}^{K+1}$ such that $y^{n}\rightarrow y$ as
$n\rightarrow
\infty$. Let $x^{n}=(x_{1}^{n},x_{2}^{n})=\Upsilon(y^{n})$ and $x=(x_{1}
,x_{2})=\Upsilon(y)$. Note that since $x\in\mathbb{D}^{K+1}$ there exists
$M>0$ such that
%
%
\begin{equation}\label{xbounded}%
\Vert\Upsilon(y)\Vert_{T}<M.
\end{equation}
Let $\Lambda$ be the set of strictly increasing functions $\lambda
\dvtx\mathbb{R}_{+}\rightarrow\mathbb{R}_{+}$ with $\lambda(0)=0$,
$\lim
_{t\rightarrow\infty}\lambda(t)=\infty$, and
\[
\gamma(\lambda)=\sup_{0\leq s<t} \biggl\vert\log\frac{\lambda
(t)-\lambda
(s)}{t-s} \biggr\vert<\infty.
\]
Since $y^{n}\rightarrow y$ as $n\rightarrow\infty$ in the
$J_{1}$-topology on $\mathbb{D}^{K+1}$, it follows from Proposition~3.5.3 of
\citet{ethkur86} that there exists a sequence $\{\lambda
^{n},n\in
\mathbb{N}\}\subset\Lambda$ such that
%
%
\begin{equation}\label{gammatozero}
\lim_{n\rightarrow\infty}\gamma(\lambda^{n})=0,%
\end{equation}
and for each $T>0$
%
%
\begin{equation}\label{convergenceyn}
{\lim_{n\rightarrow\infty}}\Vert y^{n}(\cdot)-y(\lambda^{n}(\cdot
))\Vert
_{T}=0.%
\end{equation}
For each $\lambda^{n}\in\Lambda$, $\lambda^{n}(t)$ is Lipschitz
continuous in
$t$. Hence, it is differentiable almost everywhere in $t$ with respect
to the
Lebesgue measure. Furthermore, it follows from (3.5.5) of \citet
{ethkur86} that
when $\lambda^{n}$ is differential at time $t$, its derivative ${\dot
{\lambda
}}^{n}(t)$ satisfies
%
%
\begin{equation}\label{eq:lambdaDerivatie}%
|\dot{\lambda}^{n}(t)-1|\leq\gamma(\lambda^{n}).
\end{equation}

Note that, for $i=1,2$
%
%
\begin{equation}\label{diff1}%
\int_{0}^{\lambda^{n}(t)}h_{i} ( x(s) ) \,ds=\int_{0}^{t}%
h_{i} ( x(\lambda^{n}(s)) ) \dot{\lambda}^{n}(s) \,ds.
\end{equation}
By (\ref{eq:mapPhi21}) and (\ref{diff1})
%
%
\begin{eqnarray}\label{x1eq}%
x_{1}(\lambda^{n}(t)) & = & y_{1}(\lambda^{n}(t))+\int_{0}^{\lambda^{n}(t)}
h_{1} ( x(s) ) \,ds\nonumber\\
& = & y_{1}(\lambda^{n}(t))+\int_{0}^{t}h_{1} ( x(\lambda^{n}(s)) )
\dot{\lambda}^{n}(s) \,ds\nonumber\\[-8pt]\\[-8pt]
& = & y_{1}(\lambda^{n}(t))+\int_{0}^{t}h_{1} ( x(\lambda^{n}(s)) )
\,ds\nonumber\\
&&{} -\int_{0}^{t}h_{1} ( x(\lambda^{n}(s)) ) \bigl(1-\dot{\lambda}%
^{n}(s)\bigr) \,ds.\nonumber
\end{eqnarray}
Similarly, by (\ref{eq:mapPhi22}) and (\ref{diff1})
%
%
\begin{eqnarray}\label{x2eq}%
x_{2}(\lambda^{n}(t)) & = & y_{2}(\lambda^{n}(t))+\int_{0}^{t}h_{2} (
x(\lambda^{n}(s)) ) \,ds\nonumber\\
&&{}-\int_{0}^{t}h_{2} ( x(\lambda
^{n}(s)) ) \bigl(1-\dot{\lambda}^{n}(s)\bigr) \,ds\\
&&{}+g(x_{1}(\lambda^{n}(t))).\nonumber
\end{eqnarray}
By (\ref{eq:mapPhi21}) and (\ref{x1eq})
%
%
\begin{eqnarray}\label{eq:x1}%
&& \vert x_{1}^{n}(t)-x_{1}(\lambda^{n}(t)) \vert\nonumber\\
&&\qquad \leq\vert y_{1}^{n}(t)-y_{1}(\lambda^{n}(t)) \vert+\int_{0}%
^{t} \vert h_{1}(x^{n}(s))-h_{1} ( x(\lambda^{n}(s)) )
\vert \,ds\nonumber\\
&&\qquad\quad{} +\int_{0}^{t} \vert h_{1} ( x(\lambda^{n}(s)) )
-h_{1}(0) \vert|1-\dot{\lambda}^{n}(s)| \,ds\nonumber\\[-8pt]\\[-8pt]
&&\qquad\quad{}+\int_{0}^{t} \vert
h_{1}(0) \vert|1-\dot{\lambda}^{n}(s)| \,ds\nonumber\\
&&\qquad \leq\vert y^{n}(t)-y(\lambda^{n}(t)) \vert+c\int_{0}%
^{t} \vert x^{n}(s)-x(\lambda^{n}(s)) \vert \,ds\nonumber\\
&&\qquad\quad{} +c\int_{0}^{t} \vert x(\lambda^{n}(s)) \vert|1-\dot{\lambda
}^{n}(s)| \,ds+ \vert h_{1}(0) \vert\int_{0}^{t}|1-\dot{\lambda}%
^{n}(s)| \,ds.\nonumber
\end{eqnarray}
By (\ref{eq:mapPhi22}), (\ref{x2eq}) and (\ref{eq:x1})
%
%
\begin{eqnarray}\label{eq:x2}\qquad
&& \vert x_{2}^{n}(t)-x_{2}(\lambda^{n}(t)) \vert\nonumber\\
&&\qquad \leq\vert y_{2}^{n}(t)-y_{2}(\lambda^{n}(t)) \vert\nonumber\\
&&\qquad\quad{}+\int_{0}%
^{t} \vert h_{2}(x^{n}(s))-h_{2} ( x(\lambda^{n}(s)) )
\vert \,ds \nonumber\\
&&\qquad\quad{} + \vert g(x_{1}^{n}(t))-g(x_{1}(\lambda^{n}%
(t))) \vert\nonumber\\
&&\qquad\quad{} +\int_{0}^{t} \vert h_{2} ( x(\lambda^{n}(s)) )
-h_{2}(0) \vert|1-\dot{\lambda}^{n}(s)| \,ds\nonumber\\
&&\qquad\quad{}+\int_{0}^{t} \vert
h_{2}(0) \vert|1-\dot{\lambda}^{n}(s)| \,ds\nonumber\\
&&\qquad \leq\vert y^{n}(t)-y(\lambda^{n}(t)) \vert+c\int_{0}%
^{t} \vert x^{n}(s)-x(\lambda^{n}(s)) \vert \,ds\nonumber\\
&&\qquad\quad{} +c \vert
x_{1}^{n}(t)-x_{1}(\lambda^{n}(t)) \vert+c\int_{0}^{t} \vert x(\lambda^{n}(s)) \vert|1-\dot{\lambda
}^{n}(s)| \,ds\nonumber\\
&&\qquad\quad{} + \vert h_{2}(0) \vert\int_{0}^{t}|1-\dot{\lambda}%
^{n}(s)| \,ds\\
&&\qquad \leq(1+c) \vert y^{n}(t)-y(\lambda^{n}(t)) \vert+(c+c^{2}%
)\int_{0}^{t} \vert x^{n}(s)-x(\lambda^{n}(s)) \vert \,ds\nonumber\\
&&\qquad\quad{} +(c+c^{2})\int_{0}^{t} \vert x(\lambda^{n}(s)) \vert
|1-\dot{\lambda}^{n}(s)| \,ds\nonumber\\
&&\qquad\quad{} +\bigl( \vert h_{2}(0) \vert+c \vert
h_{1}(0) \vert\bigr)\int_{0}^{t}|1-\dot{\lambda}^{n}(s)| \,ds.\nonumber
\end{eqnarray}
Then (\ref{eq:x1}) and (\ref{eq:x2}) yield
%
%
\begin{eqnarray}\label{eq:ineq}
&& \vert\Upsilon(y^{n})(t)-\Upsilon(y)(\lambda^{n}(t)) \vert
\nonumber\\
&&\qquad \leq(1+c) \vert y^{n}(t)-y(\lambda^{n}(t)) \vert\nonumber\\
&&\qquad\quad{}+(c+c^{2}%
)\int_{0}^{t} \vert\Upsilon(y^{n})(s)\,ds-\Upsilon(y)(\lambda
^{n}(s)) \vert \,ds\\
&&\qquad\quad{} +(c+c^{2})\int_{0}^{t}|1-\dot{\lambda}^{n}(s)| \vert
\Upsilon(y)(\lambda^{n}(s)) \vert \,ds\nonumber\\
&&\qquad\quad{} +\bigl( \vert h_{2}(0) \vert
+c \vert h_{1}(0) \vert\bigr)\int_{0}^{t}|1-\dot{\lambda}^{n}(s)|
\,ds.\nonumber
\end{eqnarray}
It follows from (\ref{xbounded}), (\ref{gammatozero}),
(\ref{eq:lambdaDerivatie}) and the dominated convergence theorem that
%
%
\begin{equation}\label{lastTermtoZero}%
\int_{0}^{t}|1-\dot{\lambda}^{n}(s)| \vert\Upsilon(y)(\lambda
^{n}(s)) \vert \,ds\rightarrow0 \qquad\mbox{as }n\rightarrow
\infty.
\end{equation}
Given $\delta>0$, by (\ref{gammatozero}), (\ref{eq:lambdaDerivatie})
and (\ref{lastTermtoZero}), for $n$
large enough
\begin{eqnarray*}
&&(c+c^{2})\int_{0}^{T}|1-\dot{\lambda}^{n}(s)| \vert\Upsilon
(y)(\lambda^{n}(s)) \vert \,ds\\
&&\qquad{}+\bigl( \vert h_{2}(0) \vert
+c \vert h_{1}(0) \vert\bigr)\int_{0}^{T}|1-\dot{\lambda}^{n}%
(s)| \,ds<\frac{\delta}{2}
\end{eqnarray*}
and by (\ref{convergenceyn})
\[
(1+c) \Vert y^{n}(\cdot)-y(\lambda^{n}(\cdot)) \Vert_{T}%
<\frac{\delta}{2}.
\]
By Corollary 11.2 in \citet{MandelStrong} and (\ref{eq:ineq})
\[
\Vert\Upsilon(y^{n})(\cdot)-\Upsilon(y)(\lambda^{n}(\cdot)) \Vert
_{T}\leq\delta\exp\bigl((c+c^{2})T\bigr)
\]
for large enough $n$. Thus, for each $T>0$,
\[
{\lim_{n\rightarrow\infty}} \Vert\Upsilon(y^{n})(\cdot)-\Upsilon
(y)(\lambda^{n}(\cdot)) \Vert_{T}=0.
\]
Hence, $\Upsilon(y^{n})\rightarrow\Upsilon(y)$ as $n\rightarrow
\infty$ in
$\mathbb{D}^{K+1}$ in the $J_{1}$-topology. This implies part (c)
of the lemma.

To prove part (d) of the lemma, for $y\in\mathbb{D}^{K+1}$, assume
that $x$
and $y$ satisfy (\ref{eq:mapPhi21}) and (\ref{eq:mapPhi22}). Then,
for $a>0$,
one can check that $ax$ and $ay$ also satisfy (\ref{eq:mapPhi21}) and
(\ref{eq:mapPhi22}) because of the positive homogeneity of $h_{1}$, $h_{2}$
and $g$. Therefore, $\Upsilon(ay)=a\Upsilon(y)$.\vadjust{\goodbreak}
\end{pf}

\section{\texorpdfstring{Proofs of Lemmas \protect\lowercase{\ref{thm:ssc},
\protect\ref{lem:aq}}
and Theorem \protect\lowercase{\ref{thm:virtualwaiting}}}{Proofs of Lemmas 2, 3 and Theorem 3}}
\label{sec:proof-theorems}

This section is devoted to proving Lemmas \ref{thm:ssc}, \ref{lem:aq} and
Theorem \ref{thm:virtualwaiting}. We first present two lemmas.

The first lemma is an immediate result by Proposition 4.4 of \citet
{DaiHe09}.
It proves that the virtual waiting time processes (see the
paragraph prior to Theorem~\ref{thm:virtualwaiting} for the definition)
converge to zero in distribution as $n \rightarrow\infty$.
\begin{lemma}
\label{lem:Prop4} Under the conditions of Theorem \ref{thm:critical},
\[
W^{n}\Rightarrow0 \qquad\mbox{as }n\rightarrow\infty.
\]
\end{lemma}

For $t\geq0$, let
%
%
\begin{equation}
\label{def:zeta}
\zeta^{n}(t)=\inf\{s\geq0\dvtx s+W^{n}(s)>t\}.
\end{equation}
Since $s+W^{n}(s)\leq t$ for all $s<\zeta^{n}(t)$, each customer arriving
before time $\zeta^{n}(t)$ cannot be waiting in queue at time $t$ [see Lemmas
3.2 and 3.3 of \citet{DaiHe09} for a detailed explanation];
similarly, since
$s+W^{n}(s)>t$ for all $s>\zeta^{n}(t)$, a customer who arrives after time
$\zeta^{n}(t)$ cannot be in service at $t$. So $\zeta^{n}(t)$ is a crucial
epoch with respect to the queue length at time $t$. The next lemma
concerns the process $\zeta^{n}=\{\zeta^{n}(t),t\geq0\}$.
\begin{lemma}
\label{lem:Zeta}Under the conditions of Theorem \ref{thm:critical},
$\zeta
^{n}\in\mathbb{D}$ is nondecreasing for each $n\in\mathbb{N}$, and%
\[
\zeta^{n}\Rightarrow\zeta\qquad\mbox{as }n\rightarrow\infty,
\]
where $\zeta(t)=t$ for $t\geq0$ is the identity function on $\mathbb{R}_{+}$.
\end{lemma}
\begin{pf}
First note that
%
%
\begin{equation} \label{grtt}%
\zeta^{n}(t)+W^{n}(\zeta^{n}(t))\geq t \qquad\mbox{for }t\geq0,
\end{equation}
because $W^{n}$ is right-continuous.

Next, we prove that $\zeta^{n}$ is nondecreasing in $t$. Suppose, on the
contrary, that for some $0\leq s<t$, we have $\zeta^{n}(t)<\zeta
^{n}(s)$. This
implies by (\ref{def:zeta}) that for any $\zeta^{n}(t)<u<\zeta^{n}(s)$,
\[
t<u+W^{n}(u)\leq s,
\]
leading to a contradiction.

Now we prove that $\zeta^{n}\in\mathbb{D}$, that is, $\zeta^{n}$ is
right-continuous on $[0,\infty)$ and has left limits on $(0,\infty)$. Since
$\zeta^{n}(t)\leq t$ by (\ref{def:zeta}) and $\zeta^{n}$ is nondecreasing,
$\zeta^{n}(t-)$ exists for each $t>0$; therefore, $\zeta^{n}$ has
left limits
on $(0,\infty)$. To prove right-continuity, fix $\varepsilon>0$ and
$t\geq0$.
We have
\[
\zeta^{n}(t)+\varepsilon+W^{n}\bigl(\zeta^{n}(t)+\varepsilon\bigr)>t+\delta\qquad
\mbox{for some }\delta>0,
\]
so that $\zeta^{n}(t+\delta^{\prime})\leq\zeta^{n}(t+\delta)\leq
\zeta
^{n}(t)+\varepsilon$ for $0<\delta^{\prime}\leq\delta$. Hence,
$\zeta
^{n}$ is
right-continuous at $t$, proving $\zeta^{n}\in\mathbb{D}$.

Finally, we prove the convergence. By (\ref{grtt}) and the fact $\zeta
^{n}(t)\leq t$, for any $T>0$,%
\[
{\sup_{0\leq t\leq T}} \vert t-\zeta^{n}(t) \vert\leq\sup_{0\leq
t\leq
T}W^{n}(\zeta^{n}(t))\leq\sup_{0\leq t\leq T}W^{n}(t).
\]
Then $\zeta^{n}\Rightarrow\zeta$ as $n\rightarrow\infty$ follows from
Lemma \ref{lem:Prop4}.
\end{pf}
\begin{pf*}{Proof of Lemma \ref{thm:ssc}} Fix $T>0$, and restrict
$t\in
\lbrack0,T]$. Since each
customer arriving before time $\zeta^{n}(t)$ will either have entered service
or abandoned the system by time $t$, we have $(X^{n}(t))^{+}\leq
E^{n}(t)-E^{n}(\zeta^{n}(t))+\Delta^{n}(\zeta^{n}(t))$ where $\Delta
^{n}(t)=E^{n}(t)-E^{n}(t-)$ is the number of customers who arrive
(exactly) at
time $t$. Because $\zeta^{n}(t)\leq t$ by (\ref{def:zeta}), we have
\[
\sup_{0\leq t\leq T}\Delta^{n}(\zeta^{n}(t))\leq\Vert\Delta
^{n}\Vert_{T}
\]
for $\Vert\Delta^{n}\Vert_{T}=\sup_{0\leq t\leq T}\Delta^{n}(t)$; thus,
%
%
\begin{equation}\label{13}%
(X^{n}(t))^{+}\leq E^{n}(t)-E^{n}(\zeta^{n}(t))+\Vert\Delta^{n}\Vert_{T}.
\end{equation}
Similarly, because a customer who arrives during $(\zeta^{n}(t),t]$ will
either be waiting in queue at time $t$ or has abandoned the system by
$t$, one
has%
%
%
\begin{equation}\label{14}%
(X^{n}(t))^{+}\geq E^{n}(t)-E^{n}(\zeta^{n}(t))-\bigl(A^{n}(t)-A^{n}(\zeta
^{n}(t))\bigr).
\end{equation}

Let $\{\psi^{0}(i),i\in\mathbb{N}\}$ be a sequence of i.i.d.
$K$-dimensional random vectors such that for $k=1,\ldots,K$, the
probability that $\psi^{0}(i)=e^{k}$ is $p_{k}$; it is used to
indicate the initial service phase of each customer (see the first
paragraph in Section \ref{sec:dynamics}). Write
\[
\Psi^{0}(N)=\sum_{i=1}^{N}\psi^{0}(i) \quad\mbox{and}\quad \hat{\Psi}%
^{0}(N)=\Psi^{0}(N)-pN.
\]
Because the customers who arrive before time $\zeta^{n}(t)$ cannot be waiting
in queue at time $t$ (they have either abandoned the system or started
service), for $k=1,\ldots,K$,%
%
%
\begin{eqnarray} \label{15}\hspace*{20pt}
Q_{k}^{n}(t) & \leq &\Psi_{k}^{0} \bigl( (X^{n}(0))^{+}+E^{n}(t) \bigr)
-\Psi_{k}^{0} \bigl( (X^{n}(0))^{+}+E^{n}(\zeta^{n}(t))-\Vert\Delta^{n}%
\Vert_{T} \bigr) \nonumber\\
&=&\hat{\Psi}_{k}^{0} \bigl( (X^{n}(0))^{+}+E^{n}(t) \bigr) -\hat{\Psi}%
_{k}^{0} \bigl( (X^{n}(0))^{+}+E^{n}(\zeta^{n}(t))-\Vert\Delta^{n}\Vert
_{T} \bigr) \\
&&{} +p_{k} \bigl( E^{n}(t)-E^{n}(\zeta^{n}(t))+\Vert\Delta^{n}\Vert
_{T} \bigr) .\nonumber
\end{eqnarray}
Similarly, the customers who arrive during $(\zeta^{n}(t),t]$ cannot
get into
service by time $t$. Then%
%
%
\begin{eqnarray} \label{16}%
&&
Q_{k}^{n}(t)+\bigl(A^{n}(t)-A^{n}(\zeta^{n}(t))\bigr) \nonumber\\
&&\qquad \geq\Psi_{k}^{0} \bigl(
(X^{n}(0))^{+}+E^{n}(t) \bigr) -\Psi_{k}^{0} \bigl( (X^{n}(0))^{+}%
+E^{n}(\zeta^{n}(t)) \bigr) \nonumber\\[-8pt]\\[-8pt]
&&\qquad =\hat{\Psi}_{k}^{0} \bigl( (X^{n}(0))^{+}+E^{n}(t) \bigr) -\hat{\Psi}%
_{k}^{0} \bigl( (X^{n}(0))^{+}+E^{n}(\zeta^{n}(t)) \bigr) \nonumber\\
&&\qquad\quad{} +p_{k} \bigl( E^{n}(t)-E^{n}(\zeta^{n}(t)) \bigr)
.\nonumber
\end{eqnarray}
Combining (\ref{13})--(\ref{16}), we have%
%
%
\begin{equation}\label{17}%
\Lambda_{k}^{n}(t)\leq Q_{k}^{n}(t)-p_{k}(X^{n}(t))^{+}\leq\Pi_{k}^{n}(t),
\end{equation}
where%
\begin{eqnarray*}
\Lambda_{k}^{n}(t) & = & \hat{\Psi}_{k}^{0} \bigl( (X^{n}(0))^{+}%
+E^{n}(t) \bigr) -\hat{\Psi}_{k}^{0} \bigl( (X^{n}(0))^{+}+E^{n}(\zeta
^{n}(t)) \bigr) \\
&&{} -\bigl(A^{n}(t)-A^{n}(\zeta^{n}(t))\bigr)-p_{k}\Vert\Delta^{n}\Vert_{T},\\
\Pi_{k}^{n}(t) & = & \hat{\Psi}_{k}^{0} \bigl( (X^{n}(0))^{+}+E^{n}(t) \bigr)
-\hat{\Psi}_{k}^{0} \bigl( (X^{n}(0))^{+}+E^{n}(\zeta^{n}(t))-\Vert
\Delta
^{n}\Vert_{T} \bigr) \\
&&{} +p_{k} \bigl( \Vert\Delta^{n}\Vert_{T}+A^{n}(t)-A^{n}(\zeta
^{n}(t)) \bigr) .
\end{eqnarray*}
Let $\tilde{\Psi}^{0,n}(t)=\hat{\Psi}^{0}(\lfloor nt\rfloor)/\sqrt{n}$,
$\Vert\tilde{\Delta}^{n}\Vert_{T}=\Vert\Delta^{n}\Vert_{T}/\sqrt
{n}$, and
$\Vert\bar{\Delta}^{n}\Vert_{T}=\Vert\Delta^{n}\Vert_{T}/n$. Rewriting
(\ref{17}) using diffusion scaling one has
%
%
\begin{equation} \label{21}%
\tilde{\Lambda}_{k}^{n}(t)\leq\frac{1}{\sqrt{n}} \bigl( Q_{k}^{n}%
(t)-p_{k}(X^{n}(t))^{+} \bigr) \leq\tilde{\Pi}_{k}^{n}(t),
\end{equation}
where%
\begin{eqnarray*}
\tilde{\Lambda}_{k}^{n}(t) & = & \tilde{\Psi}_{k}^{0,n} \bigl( (\bar{X}%
^{n}(0))^{+}+\bar{E}^{n}(t) \bigr) -\tilde{\Psi}_{k}^{0,n} \bigl( (\bar
{X}^{n}(0))^{+}+\bar{E}^{n}(\zeta^{n}(t)) \bigr) \\
&&{} -\bigl(\tilde{A}^{n}(t)-\tilde{A}^{n}(\zeta^{n}(t))\bigr)-p_{k}\Vert
\tilde{\Delta}^{n}\Vert_{T},\\
\tilde{\Pi}_{k}^{n}(t) & = & \tilde{\Psi}_{k}^{0,n} \bigl( (\bar{X}^{n}%
(0))^{+}+\bar{E}^{n}(t) \bigr) -\tilde{\Psi}_{k}^{0,n} \bigl( (\bar{X}%
^{n}(0))^{+}+\bar{E}^{n}(\zeta^{n}(t))-\Vert\bar{\Delta}^{n}\Vert
_{T} \bigr)
\\
&&{} +p_{k}\bigl(\Vert\tilde{\Delta}^{n}\Vert_{T}+\tilde{A}^{n}(t)-\tilde
{A}%
^{n}(\zeta^{n}(t))\bigr).
\end{eqnarray*}

Next, we show that $\tilde{\Lambda}_{k}^{n}\Rightarrow0$ and $\tilde
{\Pi
}%
_{k}^{n}\Rightarrow0$ as $n\rightarrow\infty$, which, together with
(\ref{21}), will lead to (\ref{eq:SSC}). Using (\ref{eq:Eclt}), we have
%
%
\begin{equation} \label{12}%
\Vert\tilde{\Delta}^{n}\Vert_{T}\Rightarrow0 \qquad\mbox{as
}n\rightarrow
\infty.
\end{equation}
Lemma \ref{lem:Zeta}, (\ref{eq:AbConverge}), Theorem 3.9 in
\citet
{bil99} and
the random-time-change theorem [see the lemma on page 151 of
\citet{bil99}]
yield%
%
%
\begin{equation} \label{22}%
{\sup_{0\leq t\leq T}} \vert\tilde{A}^{n}(t)-\tilde{A}^{n}(\zeta
^{n}(t)) \vert\Rightarrow0 \qquad\mbox{as }n\rightarrow\infty.
\end{equation}
By Theorem \ref{thm:fluid}, Lemma \ref{lem:Zeta} and the random-time-change
theorem,
%
%
\begin{equation}\label{20}
(\bar{X}^{n}(0))^{+}\Rightarrow0 \quad\mbox{and}\quad \bar{E}^{n}(\zeta
^{n}(\cdot))\Rightarrow\bar{E} \qquad\mbox{as }n\rightarrow\infty.
\end{equation}
Since $\tilde{\Psi}^{0,n}\Rightarrow\tilde{\Psi}^{0}$ where
$\tilde{\Psi}^{0}$
is a $K$-dimensional Brownian motion, by Theorem \ref{thm:fluid},
(\ref{12}),
(\ref{20}) and the random-time-change theorem
%
%
\begin{eqnarray}
\label{23}
&& \sup_{0\leq t\leq T} \bigl\vert\tilde{\Psi}_{k}^{0,n} \bigl( (\bar{X}%
^{n}(0))^{+}+\bar{E}^{n}(t) \bigr)\nonumber\\[-8pt]\\[-8pt]
&&\qquad\hspace*{4.3pt}{} -\tilde{\Psi}_{k}^{0,n} \bigl( (\bar
{X}^{n}(0))^{+}+\bar{E}^{n}(\zeta^{n}(t))-\Vert\bar{\Delta
}^{n}\Vert
_{T} \bigr) \bigr\vert\Rightarrow0,\nonumber\\
\label{24}%
&& \sup_{0\leq t\leq T} \bigl\vert\tilde{\Psi}_{k}^{0,n} \bigl( (\bar{X}%
^{n}(0))^{+}+\bar{E}^{n}(t) \bigr)\nonumber\\[-8pt]\\[-8pt]
&&\qquad\hspace*{4.3pt}{} -\tilde{\Psi}_{k}^{0,n} \bigl( (\bar
{X}^{n}(0))^{+}+\bar{E}^{n}(\zeta^{n}(t)) \bigr) \bigr\vert\Rightarrow0\nonumber
\end{eqnarray}
as $n\rightarrow\infty$. We deduce from (\ref{12})--(\ref{24}) that
$\tilde{\Lambda}_{k}^{n}\Rightarrow0$ and $\tilde{\Pi
}_{k}^{n}\Rightarrow0$ as
$n\rightarrow\infty$.
\end{pf*}
\begin{pf*}{Proof of Theorem \ref{thm:virtualwaiting}} Since all
customers arriving prior to time
$t\geq0$ will have either got into service or abandoned the system by time
$t+W^{n}(t)$ [see Lemmas 3.2 and 3.3 of \citet{DaiHe09}], then%
\[
\bigl(X^{n}\bigl(t+W^{n}(t)\bigr)\bigr)^{+}\leq E^{n}\bigl(t+W^{n}(t)\bigr)-E^{n}(t).
\]
For a customer who arrives during $(t,t+W^{n}(t)]$, he can possibly be waiting
in queue at time $t+W^{n}(t)$, or have abandoned the system by $t+W^{n}(t)$,
or starts his service (exactly) at $t+W^{n}(t)$. Therefore,%
\begin{eqnarray*}
E^{n}\bigl(t+W^{n}(t)\bigr)-E^{n}(t)&\leq&\bigl(X^{n}\bigl(t+W^{n}(t)\bigr)\bigr)^{+}+A^{n}\bigl(t+W^{n}%
(t)\bigr)-A^{n}(t)\\
&&{}+\Delta_{D}^{n}\bigl(t+W^{n}(t)\bigr),
\end{eqnarray*}
where $\Delta_{D}^{n}(t)=D^{n}(t)-D^{n}(t-)$ is the number of service
completions (exactly) at time $t$. Then by (\ref{eq:tilEn}) and
(\ref{eq:tilAn}),
\begin{eqnarray*}
0 & \leq &\frac{1}{\sqrt{n}}\lambda^{n}W^{n}(t)-\bigl(\tilde
{X}^{n}\bigl(t+W^{n}%
(t)\bigr)\bigr)^{+}+\tilde{E}^{n}\bigl(t+W^{n}(t)\bigr)-\tilde{E}^{n}(t)\\
& \leq& \tilde{A}^{n}\bigl(t+W^{n}(t)\bigr)-\tilde{A}^{n}(t)
+\tilde{\Delta}_{D}
^{n}\bigl(t+W^{n}(t)\bigr),
\end{eqnarray*}
where $\tilde{\Delta}_{D}^{n}(t)=\Delta_{D}^{n}(t)/\sqrt{n}$. This
leads to
%
%
\begin{eqnarray} \label{4}%
&&
\bigl\vert\mu\sqrt{n}W^{n}(t)-\bigl(\tilde{X}^{n}\bigl(t+W^{n}(t)\bigr)\bigr)^{+} \bigr\vert
\nonumber\\
&&\qquad
\leq\biggl\vert\sqrt{n} \biggl( \frac{1}{n}\lambda^{n}-\mu\biggr)%
W^{n}(t) \biggr\vert+ \bigl\vert\tilde{E}^{n}\bigl(t+W^{n}(t)\bigr)-\tilde{E}%
^{n}(t) \bigr\vert\\
&&\qquad\quad{} + \bigl\vert\tilde{A}^{n}\bigl(t+W^{n}(t)\bigr)-\tilde{A}^{n}(t) \bigr\vert
+\tilde{\Delta}_{D}^{n}\bigl(t+W^{n}(t)\bigr).\nonumber
\end{eqnarray}

Next we show that all terms on the right-hand side of (\ref{4})
converge weakly to
zero as $n\rightarrow\infty$. Using (\ref{eq:lambdan}) and
Lemma \ref{lem:Prop4}, we get%
%
%
\begin{equation} \label{7}%
\biggl\vert\sqrt{n} \biggl( \frac{1}{n}\lambda^{n}-\mu\biggr) W^{n} \biggr\vert
\Rightarrow0 \qquad\mbox{as }n\rightarrow\infty.
\end{equation}
For any $T>0$, by (\ref{eq:Eclt}) and Lemma \ref{lem:Prop4},
%
%
\begin{equation} \label{8}%
\sup_{0\leq t\leq T} \bigl\vert\tilde{E}^{n}\bigl(t+W^{n}(t)\bigr)-\tilde{E}%
^{n}(t) \bigr\vert\Rightarrow0 \qquad\mbox{as }n\rightarrow\infty.
\end{equation}
By (\ref{eq:AbConverge}) and Lemma \ref{lem:Prop4},
%
%
\begin{equation} \label{9}%
\sup_{0\leq t\leq T} \bigl\vert\tilde{A}^{n}\bigl(t+W^{n}(t)\bigr)-\tilde{A}%
^{n}(t) \bigr\vert\Rightarrow0 \qquad\mbox{as }n\rightarrow\infty.
\end{equation}
Set $\tilde{D}^{n}(t)=(D^{n}(t)-n\mu t)/\sqrt{n}$. It follows from
(\ref{eq:DnewRep}) that $\tilde{D}^{n}\Rightarrow\tilde{D}$ as
$n\rightarrow
\infty$, where
\[
\tilde{D}(t)=-e^{\prime}\tilde{M}(t)+e^{\prime}R\int_{0}^{t}\tilde
{Z}(s) \,ds.
\]
Since $\tilde{D}$ is continuous almost surely, using Lemma \ref
{lem:Prop4} again, we have%
%
%
\begin{equation} \label{10}%
\sup_{0\leq t\leq T} \bigl\vert\tilde{\Delta}_{D}^{n}\bigl(t+W^{n}(t)\bigr) \bigr\vert
\Rightarrow0 \qquad\mbox{as }n\rightarrow\infty.
\end{equation}
Combining (\ref{4})--(\ref{10}), we deduce that%
%
%
\begin{equation}\label{18}%
\sup_{0\leq t\leq T} \bigl\vert\mu\sqrt{n}W^{n}(t)-\bigl(\tilde{X}^{n}%
\bigl(t+W^{n}(t)\bigr)\bigr)^{+} \bigr\vert\Rightarrow0 \qquad\mbox{as }n\rightarrow\infty.
\end{equation}

By (\ref{eq:mapPhi1}), the process $\tilde{X}$ is continuous almost surely;
then so is $(\tilde{X})^{+}$. Because $s+W^{n}(s)\leq t+W^{n}(t)$ for
$0\leq
s\leq t$ [see Lemma 3.3 of \citet{DaiHe09}] and the process
$(\tilde
{X})^{+}$ is
continuous almost surely, by Lemma \ref{lem:Prop4} and the random-time-change
theorem,%
%
%
\begin{equation} \label{11}%
\bigl(\tilde{X}^{n}\bigl(\cdot+W^{n}(\cdot)\bigr)\bigr)^{+}\Rightarrow(\tilde{X})^{+}
\qquad\mbox{as }n\rightarrow\infty.
\end{equation}
By (\ref{18}), (\ref{11}) and the convergence-together theorem [see
Theorem 3.1
of \citet{bil99}], $\sqrt{n}W^{n}\Rightarrow(\tilde
{X})^{+}/\mu$
as $n\rightarrow\infty$.
\end{pf*}
\begin{pf*}{Proof of Lemma \ref{lem:aq}}
Recall that any customer who is waiting in queue at time $t\geq0$ must arrive
at the system during $[\zeta^{n}(t),t]$ [see (\ref{def:zeta}) and the
discussion therein], and must leave the queue (either goes into service or
abandons the system) by time $t+W^{n}(t)$ [see Lemmas 3.2 and 3.3 of
\citet{DaiHe09}]. This implies
\[
A_{Q}^{n}(t)\leq A^{n}\bigl(t+W^{n}(t)\bigr)-A^{n}(\zeta^{n}(t)-).
\]
It follows that for any $T>0$,%
\begin{eqnarray*}
\sup_{0\leq t\leq T}\tilde{A}_{Q}^{n}(t)&\leq&\sup_{0\leq t\leq T}
\bigl\vert
\tilde{A}^{n}\bigl(t+W^{n}(t)\bigr)-\tilde{A}^{n}(\zeta^{n}(t)) \bigr\vert\\
&&{}+{\sup
_{0\leq
t\leq T} }\vert\tilde{A}^{n}(\zeta^{n}(t))-\tilde{A}^{n}(\zeta
^{n}(t)-) \vert.
\end{eqnarray*}
By (\ref{22}) and (\ref{9})
\[
\sup_{0\leq t\leq T} \bigl\vert\tilde{A}^{n}\bigl(t+W^{n}(t)\bigr)-\tilde{A}^{n}%
(\zeta^{n}(t)) \bigr\vert\Rightarrow0 \qquad\mbox{as }n\rightarrow\infty.
\]
By (\ref{eq:AbConverge}) and the fact $\zeta^{n}(t)\leq t$
\begin{eqnarray*}
&&
{\sup_{0\leq t\leq T}} \vert\tilde{A}^{n}(\zeta^{n}(t))-\tilde{A}^{n}
(\zeta^{n}(t)-) \vert\\
&&\qquad\leq{\sup_{0\leq t\leq T}} \vert\tilde{A}%
^{n}(t)-\tilde{A}^{n}(t-) \vert\Rightarrow0 \qquad\mbox{as }%
n\rightarrow\infty.
\end{eqnarray*}
Hence, (\ref{eq:AQ}) holds.
\end{pf*}
\end{appendix}


%
\printaddresses

\end{document}